\documentclass[12pt]{amsart}

\usepackage{amssymb}

\usepackage{enumitem}

\usepackage{graphicx}

\makeatletter
\@namedef{subjclassname@2020}{%
  \textup{2020} Mathematics Subject Classification}
\makeatother

\usepackage[T1]{fontenc}

\newtheorem{theorem}{Theorem}[section]

\newtheorem{lemma}[theorem]{Lemma}
\newtheorem{proposition}[theorem]{Proposition}

\theoremstyle{definition}

\newtheorem{remark}[theorem]{Remark}
\newtheorem{example}[theorem]{Example}

\numberwithin{equation}{section}

\begin{document}


\baselineskip=17pt


\title[Hypersurfaces satisfying a particular Roter type equation]{Hypersurfaces in spaces of constant curvature satisfying a particular Roter type equation}

\author[R. Deszcz]{Ryszard Deszcz}
\address{retired employee of the Department of Applied Mathematics\\ 
Wroc{\l}aw University of Environmental and Life Sciences\\
Grunwaldzka 53\\
50-357 Wroc{\l}aw, Poland}
\email{Ryszard.Deszcz@upwr.edu.pl}

\author[M. G{\l}ogowska]{Ma{\l}gorzata G{\l}ogowska}
\address{Department of Applied Mathematics\\ 
Wroc{\l}aw University of Environmental and Life Sciences\\
Grunwaldzka 53\\
50-357 Wroc{\l}aw, Poland}
\email{Malgorzata.Glogowska@upwr.edu.pl}

\author[M. Hotlo\'{s}]{Marian Hotlo\'{s}}
\address{retired employee of the Department of Applied Mathematics\\ 
Wroc{\l}aw University of Science and Technology\\
Wybrze\.{z}e Wyspia\'{n}skiego 27\\ 
50-370 Wroc{\l}aw, Poland} 
\email{Marian.Hotlos@pwr.edu.pl}

\author[K. Sawicz]{Katarzyna Sawicz}
\address{Department of Mathematics\\ 
Silesian University of Technology\\
Kaszubska 23\\ 
44-100 Gliwice, Poland}
\email{Katarzyna.Sawicz@polsl.pl}


\begin{abstract}
Let $M$ be a hypersurface isometrically immersed 
in an $(n+1)$-dimensional semi-Riemannian space of constant curvature, $n > 3$, 
such that its shape operator ${\mathcal A}$ satisfies
${\mathcal A}^{3} = \phi {\mathcal A}^{2} + \psi {\mathcal A} + \rho I\!d$,
where $\phi$, $\psi$ and $\rho$ are some functions on $M$
and $I\!d$ is the identity operator.
The main result of this paper states
that on the set $U$ of all points of $M$ at which 
the square $\mathcal{S}^{2}$ of the Ricci operator $\mathcal{S}$  of $M$ is not 
a linear combination of $\mathcal{S}$ and $I\!d$,
the Riemann-Christoffel curvature tensor $R$ of $M$
is a linear combination of some Kulkarni-Nomizu products 
formed by the metric tensor $g$, the Ricci tensor $S$ 
and the tensor $S^{2}$ of $M$,
i.e., the tensor $R$ satisfies on $U$ some Roter type equation.
Moreover, the $(0,4)$-tensor $R \cdot S$ 
is on $U$ a linear combination of some Tachibana tensors formed by the tensors $g$, $S$ and $S^{2}$. 
In particular, if $M$ is a hypersurface isometrically immersed in the $(n+1)$-dimensional 
Riemannian space of constant curvature, 
$n > 3$, with three distinct principal curvatures 
and the Ricci operator $\mathcal{S}$ with three distinct eigenvalues
then the Riemann-Christoffel curvature tensor $R$ of $M$ also satisfies 
a Roter type equation of this kind.
\end{abstract}

\subjclass[2020]{Primary 53B20, 53B25}

\keywords{Einstein space, partially Einstein space,
quasi-Einstein space, 2-quasi-Einstein space,
Ricci principal curvature,
warped product manifold,
pseudosymmetry type curvature condition, Roter equation, Roter space,
Roter type equation, generalized Roter space,
hypersurface, principal curvature, generalized Clifford torus, 
type number two hypersurface, 2-quasi umbilical hypersurface, 
Cartan hypersurface, austere hypersurface.}

\maketitle

\dedicatory{Dedicated to the memory of Professor Jacek Gancarzewicz
(1945-2013)}


\section{Introduction}
A semi-Riemannian manifold $(M,g)$, $n = \dim M \geq 2$, is said to be 
an {\sl Einstein manifold} \cite{Besse-1987},
or an {\sl Einstein space}, if at every point of $M$ 
its Ricci tensor $S$ is proportional to $g$, 
i.e., $S\, =\, ( \kappa /n)\, g$ on $M$,
assumed that $\kappa$ is constant when $n = 2$. 
According to {\cite[p. 432] {Besse-1987}}
this condition is called the {\sl Einstein metric condition}.
Einstein manifolds form a natural subclass
of several classes of semi-Riemannian manifolds which are determined by
curvature conditions imposed on their Ricci tensor 
{\cite[Table, pp. 432-433] {Besse-1987}}.
These conditions are called 
{\sl generalized Einstein metric conditions}
{\cite[Chapter XVI] {Besse-1987}}.

Let $(M,g)$, $n \geq 4$, be a semi-Riemannian manifold. 
We define the following subsets of $M$:
${\mathcal U}_{R} 
= \{x \in M|\, R \neq  \frac{\kappa }{2(n-1) n } g \wedge g 
\ \mbox {at}\ x \}$, 
${\mathcal U}_{S} = \{x \in M|\, S \neq \frac{\kappa }{n}\, g
\ \mbox {at}\ x \}$ and 
${\mathcal U}_{C} = \{x \in M|\, C \neq 0\ \mbox {at}\ x \}$. 
We note that  ${\mathcal U}_{S} \cup  {\mathcal U}_{C} 
=  {\mathcal U}_{R}$ 
(see, e.g., {\cite[p. 379] {DGHP-TZ 1}}).
We refer to Section 2 of this paper for precise definitions 
of the symbols used
(see also \cite{{DGHP-TZ 1}, {DGHP-TZ 2}, {DGHSaw}, {R99}, {2018_DGZ}, {DHV2008}}).

The semi-Riemannian manifold $(M,g)$ is said to be 
a {\sl quasi-Einstein manifold} (see, e.g., \cite{DGHS}), 
or a {\sl quasi-Einstein space}, if 
\begin{eqnarray}
\mathrm{rank}(S - \alpha\, g) &=& 1
\label{quasi02}
\end{eqnarray}
on ${\mathcal U}_{S} \subset M$, where $\alpha $
is some function on ${\mathcal U}_{S}$.
It is known that every non-Einstein warped product manifold 
$\overline{M} \times _{F} \widetilde{N}$
with a $1$-dimensional base manifold $(\overline{M}, \overline{g})$ and
a $2$-dimensional manifold $(\widetilde{N}, \widetilde{g})$
or an $(n-1)$-dimensional Einstein manifold
$(\widetilde{N}, \widetilde{g})$, $n \geq 4$, and a warping function $F$,
is a  quasi-Einstein manifold (see, e.g., \cite{{Ch-DDGP}, {2016_DGJZ}}). 
Evidently, (\ref{quasi02}) implies
{\cite[eq. (8)] {G6}}
\begin{eqnarray}
S^{2} &=& ( \kappa - (n-2) \alpha )\, S 
+ \alpha ( (n-1) \alpha - \kappa  )\, g .
\label{quasi0202}
\end{eqnarray}
Quasi-Einstein manifolds arose during the study of exact solutions
of the Einstein field equations and the investigation 
on quasi-umbilical hypersurfaces of conformally flat spaces, 
see, e.g., \cite{{DGHSaw}, {2016_DGJZ}} and references therein. 
Quasi-Einstein hypersurfaces 
in semi-Riemannian space forms
were studied among others things in
\cite{{DGHS}, {2016_DGHZhyper}, {2018_DGZ}, {DHS105}, {G6}}, 
see also \cite{DGHSaw} and references therein.
Quasi-Einstein manifolds satisfying some pseudo\-sym\-met\-ry type 
curvature conditions were investigated for instance in
\cite{{Ch-DDGP}, {2015_DGHZ}, {2016_DGJZ}, {SDHJK}}.

The semi-Riemannian manifold $(M,g)$, $n \geq 3$, will be called
a {\sl{partially Einstein manifold}}, or
a {\sl{partially Einstein space}}
(cf. {\cite[Foreword] {CHEN-2017}}, 
{\cite[p. 20] {V2}}, \cite{LV3-Foreword}),  
if on ${\mathcal{U}}_{S} \subset M$ 
its Ricci operator ${\mathcal{S}}$ and
the identity operator $I\!d$ on the Lie algebra $\frak{X} (M)$ of the vector fields on $M$
satisfy
${\mathcal{S}}^{2} = \lambda {\mathcal{S}} + \mu  I\!d_{x}$,
or equivalently, 
\begin{eqnarray}
S^{2} &=&  \lambda \, S + \mu \, g ,
\label{partiallyEinstein}
\end{eqnarray}
where $\lambda$ and $ \mu$ are some functions on ${\mathcal{U}}_{S}$.
Evidently, (\ref{quasi0202}) 
is a special case of (\ref{partiallyEinstein}). 
Thus every quasi-Einstein manifold is partially Einstein.
The converse statement is not true. 
Contracting (\ref{partiallyEinstein}) we get
$ \mathrm{tr}_{g} (S^{2}) = \lambda \, \kappa + n\, \mu$, which 
together with 
(\ref{partiallyEinstein}) yields (cf. {\cite[Section 5] {2021-DGH}})
\begin{eqnarray} 
S^{2} - \frac{  \mathrm{tr}_{g} (S^{2})} {n} \, g 
&=&  \lambda \left( S - \frac{\kappa }{n} \, g \right) .
\label{partiallyEinstein.11} 
\end{eqnarray}
Thus in particular, 
if a Riemannian manifold $(M,g)$, $n \geq 3$, is a partially Einstein space
then at every point $x \in {\mathcal{U}}_{S} \subset M$ 
its Ricci operator ${\mathcal{S}}$
has exactly two distinct eigenvalues,
i.e., exactly two distinct Ricci principal curvatures,
$\kappa _{1}$ and $\kappa _{2}$ 
with multiplicities $p$ and $n-p$, respectively, where
$1 \leq p \leq n-1$.
Evidently, if $p = 1$, or $p = n-1$, 
then $(M,g)$ is a quasi-Einstein manifold.

If the Riemann-Christoffel curvature tensor $R$,
or equivalently, the Weyl conformal curvature tensor $C$,
of a non-partially Einstein and non-conformally flat 
semi-Rieman\-nian manifold  $(M,g)$, $n \geq 4$,
is expressed at every point of  
${\mathcal U}_{S} \cap {\mathcal U}_{C} \subset M$ 
by linear combination of the Kulkarni-Nomizu products: 
$g \wedge g$, $g \wedge S$ and $S \wedge S$, then such a manifold
is called
a {\sl Roter type manifold}, or a {\sl Roter manifold}, 
or a {\sl Roter space} {\cite[Section 15] {CHEN-2021}},
\cite{{DGHP-TZ 1}, {DGHP-TZ 2},
{DGHSaw}, {DGP-TV02}, {2018_DGZ}, {DHV2008}, {2018_DH}}.

The curvature tensor $R$ of a Roter space $(M,g)$, $n \geq 4$, 
satisfies on ${\mathcal U}_{S} \cap {\mathcal U}_{C} \subset M$ 
\begin{eqnarray} 
R &=& \frac{\phi}{2}\, S\wedge S + \mu\, g\wedge S 
+ \frac{\eta}{2}\, g \wedge g ,
\label{eq:h7a}
\end{eqnarray}
where $\phi$, $\mu $ and $\eta $ are some functions on this set.

It is easy to check 
that at every point of ${\mathcal U}_{S} \cap {\mathcal U}_{C} \subset M$
of a Roter space $(M,g)$, $n \geq 4$, its tensor $S^{2}$ 
is a linear combination of the tensors $S$ and $g$,
i.e., (\ref{partiallyEinstein}) holds on
${\mathcal U}_{S} \cap {\mathcal U}_{C}$
(see, e.g., {\cite[Theorem 2.4] {2018_DGZ}}, 
{\cite[Theorem 2.1] {2018_DH}}).
Precisely, (\ref{partiallyEinstein.11}) 
with $\lambda = \phi ^{-1} ( (n-2) \mu + \phi \kappa -1)$
is satisfied on  
${\mathcal U}_{S} \cap {\mathcal U}_{C}$. 
It is easy to verify that every Roter space 
is a non-quasi Einstein manifold.
For further results on Roter spaces we refer to
\cite{{DecuDH}, {DecuP-TSVer}, {2021-DGH}, {2016_DGJZ},
{2011_DGP-TV}, {DK2003}, {DePlaScher},  
{Kow01}, {Kow2}, {SDHJK}}.

Let $(M,g)$, $n \geq 4$, be a non-partially Einstein 
and non-conformally flat semi-Rieman\-nian manifold.
If its Riemann-Christoffel curvature tensor $R$,
or equivalently, its Weyl conformal curvature tensor $C$,
is at every point of  
${\mathcal U}_{S} \cap {\mathcal U}_{C} \subset M$ 
a linear combination of the Kulkarni-Nomizu products 
formed by the tensors:
$S^{0} = g$ and $S^{1} = S, \ldots , S^{p-1}, S^{p}$, 
where $p$ is some natural 
number $\geq 2$, then $(M,g)$
is called 
a {\sl generalized Roter type manifold},
or a {\sl generalized Roter manifold}, 
or a {\sl generalized Roter type space}, 
or a {\sl generalized Roter space} \cite{{DGHP-TZ 2}, {SDHJK}, {2016_SK}, {2019_SK}}. 
An equation in which the tensor $R$,
or equivalently, the tensor $C$,
is expressed by
the above linear combination of the Kulkarni-Nomizu 
products is called a {\sl{Roter type equation}}. 
For instance, 
the Roter type equation for the tensor $R$, when $p = 2$,
reads  
\begin{eqnarray} 
R &=&  \frac{\phi _{2} }{2}\, S^{2} \wedge S^{2} 
+ \phi_{1}\, S \wedge S^{2}  
+ \frac{\phi}{2}\, S \wedge S
+ \mu _{1}\, g \wedge S^{2}
+ \mu \, g \wedge S
+ \frac{ \eta }{2} \, g \wedge g ,
\label{B001simply}
\end{eqnarray}
where 
$\phi $, 
$\phi _{1}$,
$\phi _{2}$,
$\mu_{1}$, $\mu$ and 
$\eta$ 
are some functions on
${\mathcal U}_{S} \cap {\mathcal U}_{C}$. 
Because $(M,g)$ is a non-partially Einstein manifold,
at least one of the functions
$\phi _{1}$, $\phi _{2}$ and $\mu_{1}$ is a non-zero function.
Manifolds (hypersurfaces) satisfying (\ref{B001simply})
were investigated among others things in:
\cite{{2016_DGJZ}, {DGP-TV02}, {Saw-2006},
{Saw-2015}, {SDHJK}, {2016_SK}, {2019_SK}}.

A semi-Riemannian manifold $(M,g)$, $n \geq 3$, 
is called a {\sl $2$-quasi-Einstein manifold},
or a {\sl $2$-quasi-Einstein space}, if 
\begin{eqnarray}
\mathrm{rank}(S - \alpha \, g ) &\leq & 2 
\label{quasi0202weak}
\end{eqnarray}
on ${\mathcal U}_{S}$ and $\mathrm{rank} (S - \alpha \, g ) = 2$
on some open non-empty subset of ${\mathcal U}_{S}$, 
where $\alpha $ is some function on ${\mathcal U}_{S}$ 
(see, e.g., {\cite[Section 6] {2016_DGHZhyper}}, 
{\cite[Section 1] {2016_DGJZ}}, {\cite[Section 2] {DGP-TV02}}).
Every non-quasi-Einstein warped product manifold 
$\overline{M} \times _{F} \widetilde{N}$
with a $2$-dimensional base manifold $(\overline{M}, \overline{g})$ 
and a $2$-dimensional manifold $(\widetilde{N}, \widetilde{g})$
or an $(n-2)$-dimensional Einstein semi-Riemannian manifold
$(\widetilde{N}, \widetilde{g})$, when $n \geq 5$, 
and a war\-ping function $F$ satisfies (\ref{quasi0202weak})
(see, e.g., {\cite[Theorem 6.1] {2016_DGJZ}}). 
Thus some exact solutions of the Einstein field equations are 
non-conformally flat $2$-quasi-Einstein manifolds.
For instance, the Reissner-Nordstr\"{o}m spacetime, as well as
the Reissner-Nordstr\"{o}m-de Sitter type spacetimes 
are such manifolds (see, e.g.,
\cite{Kow2}). 
It seems that the Reissner-Nordstr\"{o}m spacetime
is the "oldest" example of 
a non-conformally flat 
$2$-quasi-Einstein warped pro\-duct manifold
{\cite[Section 1] {2016_DGJZ}}.
We also mention 
that some $2$-quasi-umbilical hypersurfaces
in semi-Riemannian spaces of constant curvature are $2$-quasi-Einstein manifolds
(see, e.g., \cite{DGP-TV02}).

Let $N^{n+1}_{s}(c)$, $n \geq 3$, be a semi-Riemannian space of constant curvature
with signature $(s,n+1-s)$, where 
$c = \widetilde{\kappa}/(n(n+1))$ and $\widetilde{\kappa}$ 
is its scalar curvature. 
Let $(M,g)$ be a connected hypersurface isometrically immersed 
in $N_{s}^{n+1}(c)$,  
where the metric tensor $g$ of $M$ is induced 
by the metric tensor of the ambient space. 
The Gauss equation of $M$ in $N_{s}^{n+1}(c)$ reads 
(see, e.g., \cite{{2016_DGHZhyper}, {DGP-TV02}, 
{DGPSS}, {2018_DGZ}, {Saw-2015}})
\begin{eqnarray}
R_{hijk} 
&=& \varepsilon  \, ( H_{hk}H_{ij} - H_{hj}H_{ik}) 
+ c\, ( g_{hk}g_{ij} - g_{hj}g_{ik}),\ \ \ 
\varepsilon \ =\ \pm 1 ,
\label{realC5}
\end{eqnarray}  
where
$g_{hk}$, $R_{hijk}$ and $H_{hk}$ 
are the local components of the metric tensor $g$, 
the curvature tensor $R$  
and the second fundamental tensor $H$ of $M$, respectively.

Let $M$ be a hypersurface in $N_{s}^{n+1}(c)$, $n \geq 3$.
We denote by ${\mathcal U}_{H}$ the set of all points of $M$
at which the $(0,2)$-tensor $H^{2}$ is not expressed 
by a linear combination of the second fundamental tensor
$H$ and the metric tensor $g$ of $M$, or equivalently, the square of the shape operator ${\mathcal A}^{2}$
 is not expressed by a linear combination
of the shape operator ${\mathcal A}$ and the identity operator
$I\!d$
on the Lie algebra $\frak{X} (M)$ of the vector fields on  $M$.
It is known that 
${\mathcal U}_{H} \subset {\mathcal{U}}_{S} \cap {\mathcal{U}}_{C} 
\subset M$, provided that $n \geq 4$,
see, e.g., {\cite[Section 1] {2018_DGZ}}.

In Section 3 we present results on
non-Einstein and non-conformally flat hypersurfaces in  
$N_{s}^{n+1}(c)$, $n \geq 4$,
satisfying at every point
$x \in {\mathcal{U}}_{S} \setminus {\mathcal U}_{H}$
\begin{eqnarray}
H^{2} &=& \alpha \, H + \beta \, g,\ \ \
\alpha , \beta \in \mathbb{R} ,
\label{2021.10.25.a2}
\end{eqnarray}
or, equivalently,
${\mathcal A}^{2} = \alpha {\mathcal A} + \beta I\!d$.
The main result of that section (see Theorem 3.3) states that if
a hypersurface  $M$ in a Riemannian space of constant curvature
$N^{n+1}(c)$, $n \geq 4$, has at every point exactly two distinct 
principal curvatures, $\lambda_{1}$ with multiplicity $p$ and
$\lambda_{2}$ with multiplicity $n-p$, $2 \leq p \leq n-2$,
such that
\begin{eqnarray}
(p-1) \lambda_{1} + (n-p-1) \lambda_{2} \neq 0
\label{2022.09.23.cc}
\end{eqnarray} 
then ${\mathcal{U}}_{S}  \cap {\mathcal{U}}_{C} = M$
and an equation of the form (\ref{eq:h7a}) holds on $M$,
i.e., $M$ is a Roter space. In particular, 
$M$ is a non-quasi-Einstein partially Einstein space.
Some Clifford hypersurfaces form 
a family of 
non-Einstein and non-conformally flat hypersurfaces
having at every  point exactly two distinct 
principal curvatures, $\lambda_{1}$ with multiplicity $p$ and
$\lambda_{2}$ with multiplicity $n-p$, $2 \leq p \leq n-2$,
and satisfying (\ref{2022.09.23.cc}) (see Example 3.4).

In Sections 4, 5 and 6 we investigate  
hypersurfaces $M$ in $N_{s}^{n+1}(c)$, $n \geq 4$,
with the tensor $H$ satisfying on ${\mathcal U}_{H} \subset M$ 
\begin{eqnarray}
H^{3} &=& \phi \, H^{2} + \psi \, H + \rho \, g ,
\label{real702} 
\end{eqnarray}
or, equivalently,
${\mathcal A}^{3} = \phi {\mathcal A}^{2} + \psi {\mathcal A} + \rho I\!d$,
where $\phi$, $\psi$ and $\rho$ are some functions on ${\mathcal U}_{H}$
and $I\!d$ is the identity transformation on ${\mathcal U}_{H}$.
Curvature properties of pseudosymmetry type of such hypersurfaces
were investigated among others in 
\cite{{DG90}, {DGHS}, {R99}, {DGP-TV02}, {DGPSS},
{2018_DGZ}, {DHS105}, {DeVerYap}, 
{G6}, {Saw-2004}, {Saw-2005}, {Saw-2006}, {Saw-2007}, {Saw-2015}}.

In Proposition 4.1 we state 
that on the set ${\mathcal U}_{H}$
of every hypersurface $M$ in $N_{s}^{n+1}(c)$, $n \geq 4$,
satisfying (\ref{real702})
a family of curvature conditions is satisfied. Among other things
\begin{eqnarray}
\tau ^{2}\, R &=& \frac{\varepsilon}{2}\, 
( S^{2} - \rho _{1}\, S - \rho _{0}\, g) 
\wedge 
( S^{2} - \rho _{1}\, S - \rho _{0}\, g) 
+ \frac{\tau ^{2} c}{2}\, g \wedge g 
\label{2022.09.19.aa}
\end{eqnarray}
on ${\mathcal U}_{H}$,  
where $\tau$, $\rho _{0}$ and $\rho _{1}$ 
are functions defined on ${\mathcal U}_{H}$ by 
(\ref{real705}), (\ref{real707rho0}) and (\ref{real707rho1}),
respectively. Evidently, if $\tau$ is non-zero at every point 
of ${\mathcal U}_{H}$ then from (\ref{2022.09.19.aa}) 
we obtain a Roter type equation of the form (\ref{B001simply}).

Let $M$ be a hypersurface in a Riemannian space of constant curvature
$N^{n+1}(c)$, $n \geq 3$,
having at every point $x \in {\mathcal U}_{H} \subset M$ 
three distinct principal curvatures  
$\lambda _{0}$, $\lambda _{1}$ and $\lambda _{2}$,
with multiplicities $n_{0}$, $n_{1}$ and $n_{2}$, respectively. 
Moreover, let the eigenvalues (the Ricci principal curvatures)
$\kappa _{0}$, $\kappa _{1}$  and $\kappa _{2}$ of the Ricci operator 
$\mathcal{S}$ of $M$
satisfy (\ref{2021-10-18-a}) on ${\mathcal U}_{H}$.
Then the function $\tau$, defined by (\ref{real705}),
is expressed on ${\mathcal U}_{H}$ by (see Theorem 5.2)
\begin{eqnarray}
\tau &=&
\frac{(\kappa _{0} - \kappa _{1}) ( \kappa _{0} - \kappa _{2}) 
( \kappa _{1} - \kappa _{2})}{(\lambda _{0} - \lambda _{1}) 
( \lambda _{0} - \lambda _{2}) ( \lambda _{1} - \lambda _{2})} .
\label{2021-11-15-a} 
\end{eqnarray}

If $M$ is a partially Einstein hypersurface in $N_{s}^{n+1}(c)$, 
$n \geq 4$, satisfying (\ref{real702}) on ${\mathcal U}_{H} \subset M$
then $S^{2} = \rho _{1}\, S + \rho _{0}\, g$ and $\tau = 0$
hold on ${\mathcal U}_{H}$ (Theorem 5.1).

An example of a partially Einstein hypersurface $M$ in a Riemannian space
of constant curvature $N^{n+1}(c)$, $n \geq 3$, having at every point exactly three
distinct principal curvatures is presented in Example 5.4 (a class of austere hypersurfaces including the Cartan hypersurfaces) 
and Example 5.5 (some type number two hypersurfaces).
In view of Theorem 5.1, $\tau = 0$ on those hypersurfaces.
Examples of hypersurfaces $M$ in an Euclidean space
$\mathbb{E}^{n+1}(c)$, $n \geq 4$, or $n \geq 5$,
having at every point exactly three
distinct principal curvatures with non-zero
function $\tau$ are presented in Examples 6.2 (some $2$-quasi-umbilical hypersurfaces)
and 6.3 (a family of hypersurfaces with three non-zero principal curvature 
with multiplicities $1$, $n_{1} \geq 2$, $n_{2} \geq 2$). 
Thus a Roter type equation of the form (\ref{B001simply})
is satisfied on those hypersurfaces.

\section{Preliminaries.}
Throughout this paper, all manifolds are assumed 
to be connected paracompact
mani\-folds of class $C^{\infty }$. Let $(M,g)$, $n \geq 3$,
be a semi-Riemannian manifold, and let $\nabla$ 
be its Levi-Civita connection and $\mathfrak{X} (M)$ 
the Lie algebra of vector fields on $M$. 
We denote by
$I\!d$
the identity operator on $\frak{X} (M)$.
We define on $M$ the endomorphisms 
$X \wedge _{A} Y$ and ${\mathcal{R}}(X,Y)$ of $\mathfrak(M)$ by
\begin{eqnarray*}
(X \wedge _{A} Y)Z \ =\ A(Y,Z)X - A(X,Z)Y ,\ \ \ 
{\mathcal R}(X,Y)Z 
\ =\ 
\nabla _X \nabla _Y Z - \nabla _Y \nabla _X Z - \nabla _{[X,Y]}Z ,
\end{eqnarray*}
respectively, where $X, Y, Z \in \mathfrak{X} (M)$ 
and 
$A$ is a symmetric $(0,2)$-tensor on $M$. 
The Ricci tensor $S$, the Ricci operator ${\mathcal{S}}$ 
and the scalar curvature $\kappa $ of $(M,g)$ are defined by 
\begin{eqnarray*}
S(X,Y)\ =\ \mathrm{tr} \{ Z \rightarrow {\mathcal{R}}(Z,X)Y \} ,\ \ \ 
g({\mathcal S}X,Y)\ =\ S(X,Y) ,\ \ \ 
\kappa \ =\ \mathrm{tr}\, {\mathcal{S}},
\end{eqnarray*}
respectively. The endomorphism ${\mathcal{C}}(X,Y)$ is defined by
\begin{eqnarray*}
{\mathcal C}(X,Y)Z  &=& {\mathcal R}(X,Y)Z 
- \frac{1}{n-2} \left(X \wedge _{g} {\mathcal S}Y 
+ {\mathcal S}X \wedge _{g} Y 
- \frac{\kappa}{n-1}X \wedge _{g} Y \right) Z .
\end{eqnarray*}
Now the Riemann-Christoffel curvature tensor $R$ and
the Weyl conformal curvature tensor $C$ of $(M,g)$ are defined by
\begin{eqnarray*}
R(X_1,X_2,X_3,X_4) \ =\ g({\mathcal R}(X_1,X_2)X_3,X_4) ,\ \ \
C(X_1,X_2,X_3,X_4) \ =\ g({\mathcal C}(X_1,X_2)X_3,X_4) ,
\end{eqnarray*}
respectively, where $X_1,X_2,\ldots \in \mathfrak{X}(M)$.

Let ${\mathcal B}$ be a tensor field sending any $X, Y \in \mathfrak{X} (M)$ 
to a skew-symmetric endomorphism ${\mathcal B}(X,Y)$, 
and let $B$ be the $(0,4)$-tensor associated with ${\mathcal B}$ by
\begin{eqnarray}
B(X_1,X_2,X_3,X_4) &=& 
g({\mathcal B}(X_1,X_2)X_3,X_4) .
\label{DS5}
\end{eqnarray}
The tensor $B$ is said to be a {\sl{generalized curvature tensor}}  if the
following two conditions are fulfilled:
$B(X_1,X_2,X_3,X_4) = B(X_3,X_4,X_1,X_2)$ and
\begin{eqnarray*}
B(X_1,X_2,X_3,X_4) + B(X_2,X_3,X_1,X_4) + B(X_3,X_1,X_2,X_4) &=& 0 . 
\end{eqnarray*}

For ${\mathcal B}$ as above, let $B$ be again defined by (\ref{DS5}). 
We extend the endomorphism ${\mathcal B}(X,Y)$
to a derivation ${\mathcal B}(X,Y) \cdot \, $ of the algebra 
of tensor fields on $M$,
assuming that it commutes with contractions and 
${\mathcal B}(X,Y) \cdot \, f  = 0$ for any smooth function $f$ on $M$. 
Now for a $(0,k)$-tensor field $T$,
$k \geq 1$, we can define the $(0,k+2)$-tensor $B \cdot T$ by
\begin{eqnarray*}
& & (B \cdot T)(X_1,\ldots ,X_k,X,Y) \ =\ 
({\mathcal B}(X,Y) \cdot T)(X_1,\ldots ,X_k)\\  
&=& - T({\mathcal{B}}(X,Y)X_1,X_2,\ldots ,X_k)
- \cdots - T(X_1,\ldots ,X_{k-1},{\mathcal{B}}(X,Y)X_k) .
\end{eqnarray*}
If $A$ is a symmetric $(0,2)$-tensor then we define the
$(0,k+2)$-tensor $Q(A,T)$ by
\begin{eqnarray*}
& & Q(A,T)(X_1, \ldots , X_k, X,Y) \ =\
(X \wedge _{A} Y \cdot T)(X_1,\ldots ,X_k)\\  
&=&- T((X \wedge _A Y)X_1,X_2,\ldots ,X_k) 
- \cdots - T(X_1,\ldots ,X_{k-1},(X \wedge _A Y)X_k) .
\end{eqnarray*}
In this manner we obtain the $(0,6)$-tensors 
$B \cdot B$ and $Q(A,B)$. Substituting 
${\mathcal{B}} = {\mathcal{R}}$ or ${\mathcal{B}} = {\mathcal{C}}$, 
$T=R$ or $T=C$ or $T=S$, $A=g$ or $A=S$ in the above formulas, 
we get the tensors $R\cdot R$, $R\cdot C$, $C\cdot R$, $R\cdot S$, 
$Q(g,R)$, $Q(S,R)$, $Q(g,C)$ and $Q(g,S)$.

For a symmetric $(0,2)$-tensor $A$ we denote by 
${\mathcal{A}}$ the endomorphism related to $A$ by 
$g({\mathcal{A}}X,Y) = A(X,Y)$.
The tensors
$A^{p}$, $p = 2, 3, \ldots $, are defined by 
$A^{p}(X,Y) = A^{p-1} ({\mathcal{A}}X, Y)$,
assuming that $A^{1} = A$. In this way, if $A = S$ and 
${\mathcal{A}} = {\mathcal S}$,
we obtain the tensors $S^{p}$, 
$p = 2, 3, \ldots $, assuming that $S^{1} = S$.

For a symmetric $(0,2)$-tensor $E$ and a $(0,k)$-tensor $T$, $k \geq 2$, we
define their {\sl{Kulkarni-Nomizu tensor}} $E \wedge T$ by 
(see, e.g., {\cite[Section 2] {DGHP-TZ 1}}) 
\begin{eqnarray*}
& &(E \wedge T )(X_{1}, X_{2}, X_{3}, X_{4}; Y_{3}, \ldots , Y_{k})\\
&=&
E(X_{1},X_{4}) T(X_{2},X_{3}, Y_{3}, \ldots , Y_{k})
+ E(X_{2},X_{3}) T(X_{1},X_{4}, Y_{3}, \ldots , Y_{k} )\\
& &
- E(X_{1},X_{3}) T(X_{2},X_{4}, Y_{3}, \ldots , Y_{k})
- E(X_{2},X_{4}) T(X_{1},X_{3}, Y_{3}, \ldots , Y_{k}) .
\end{eqnarray*}

It is obvious that the following tensors 
are generalized curvature tensors: $R$, $C$ and 
$E \wedge F$, where $E$ and $F$ are symmetric $(0,2)$-tensors. 
We have 
\begin{eqnarray}
C &=& R - \frac{1}{n-2}\, g \wedge S + \frac{\kappa }{2 (n-2) (n-1)} \, g \wedge g , 
\label{Weyl}\\
Q(E,F) &=& - Q(F,E)
\label{2022.09.16.a}
\end{eqnarray}
and (see, e.g., {\cite[Section 2] {DGHP-TZ 1}})
\begin{eqnarray}
\ \ \
(a)
\ \
Q(E, E \wedge F) \ =\ - \frac{1}{2}\, Q(F, E \wedge E ) ,
\ \ \
(b)
\ \
E \wedge Q(E,F) \ =\ - \frac{1}{2}\, Q(F, E \wedge E ).
\label{DS7}
\end{eqnarray}
We also have the following identity
(see, e.g., 
{\cite[Lemma 2.1(i)] {Ch-DDGP}} and references therein) 
\begin{eqnarray}
E_{1} \wedge Q(E_{2},F) + E_{2} \wedge Q(E_{1},F) 
+ Q(F, E_{1} \wedge E_{2}) 
&=& 0 ,
\label{2021.11.30.a4}
\end{eqnarray} 
where $E_{1}$, $E_{2}$ and $F$
are symmetric $(0,2)$-tensors.
We note that 
(\ref{2022.09.16.a}) and (\ref{2021.11.30.a4}) yield
\begin{eqnarray*}
Q(F, E_{1} \wedge E_{2} ) + Q(E_{1}, E_{2} \wedge F ) 
+ Q(E_{2}, F \wedge E_{1} ) &=& 0 .
\end{eqnarray*}

Let $T$ be a $(0,k)$-tensor, $k = 2, 3, \ldots $.  
The tensor $Q(A,T)$ is called the {\sl{Tachibana tensor}} of $A$ and $T$, 
or the Tachibana tensor for short (see, e.g., \cite{DGPSS}). 
By an application of (\ref{DS7})(a) we obtain on $M$ the identities
\begin{eqnarray*}
Q(g, g \wedge S) &=& - \frac{1}{2}\, Q(S,g \wedge g) ,
\ \ \ Q(S, g \wedge S)\ = \ - \frac{1}{2}\, Q(g, S \wedge S) .  
\end{eqnarray*}

On every semi-Riemannian manifold $(M,g)$, $n \geq 4$,
the following identity is satisfied 
{\cite[eq. (1.13), Theorem 3.4(i)] {2016_DGJZ}}
\begin{eqnarray}
C \cdot R + R \cdot C &=& C \cdot C  + R \cdot R 
- \frac{1}{(n-2)^{2}} \, 
Q\left( g, g \wedge \left( S^{2} 
- \frac{\kappa }{n-1}\, S \right) \right) .
\label{2021.11.30.a1}
\end{eqnarray} 
From (\ref{2021.11.30.a1}), by a suitable contraction, we get 
(cf.  {\cite[Lemma 2.3] {DGP-TV02}})
\begin{eqnarray*}
C \cdot S 
&=& R \cdot S - \frac{1}{n-2} 
Q\left( g, S^{2} - \frac{\kappa }{n-1} S \right) .
\end{eqnarray*}

We mention that the following condition 
is satisfied on every hypersurface $M$ 
in $N_{s}^{n+1}(c)$, $n \geq 4$,  
(see, e.g., {\cite[eq. (14)] {DGP-TV02}}) 
\begin{eqnarray}
R \cdot R &=& Q(S,R) - (n-2)c\, Q(g,C) .
\label{2021.10.25.a1}
\end{eqnarray} 
Thus 
on every hypersurface $M$ 
in $N_{s}^{n+1}(c)$, $n \geq 4$, 
by making use of (\ref{2021.10.25.a1}), the
identity (\ref{2021.11.30.a1}) 
turns into (cf. {\cite[Theorem 3.4(ii)] {2016_DGJZ}})
\begin{eqnarray*}
C \cdot R + R \cdot C 
&=& 
C \cdot C + Q(S,R) 
- (n-2)c\, Q(g,C)\nonumber\\
& &
- \frac{1}{(n-2)^{2}} \, 
Q\left( g, g \wedge 
\left( S^{2} - \frac{\kappa }{n-1}\, S \right) \right) .
\end{eqnarray*}

Let  $(M,g)$, $n \geq 4$, be a Roter space.
We can check that at every point of 
${\mathcal U}_{S} \cap {\mathcal U}_{C} \subset M$
the tensors $R \cdot R$ and $Q(g,R)$ are lineary dependent 
(see, e.g., {\cite[Theorem 6.7] {DGHSaw}}, 
{\cite[Theorem 3.2] {2016_DGJZ}}, {\cite[Theorem 2.1] {2018_DH}}). 
Semi-Riemannian manifolds with that property are called pseudosymmetric
manifolds. More precisely,  
a semi-Riemannian manifold $(M,g)$, $n \geq 3$, is said to be
{\sl pseudosymmetric} if the tensors $R \cdot R$ and $Q(g,R)$ 
are linearly dependent at every point of $M$
(see, e.g., 
{\cite[Chapter 8.5.3] {M7}},
{\cite[Section 15] {CHEN-2021}}, 
{\cite[Chapter 12.4] {M6}}, 
\cite{{DGHSaw}, {DHV2008}, {DVV1991}, 
{SDHJK}, {V2}, {LV2}, {LV3-Foreword}, {LV-2020}}). 
This is equivalent to
\begin{eqnarray}
R \cdot R &=& L_{R}\, Q(g,R) 
\label{pseudo}
\end{eqnarray}
on ${\mathcal{U}}_{R}$,
where $L_{R}$ is some function on ${\mathcal U}_{R}$.
We note that (\ref{pseudo}) implies
\begin{eqnarray*}
(a)\ \ R \cdot S \ =\ L_{R}\, Q(g,S)\ \ &\mbox{and}& \ \ 
(b)\ \ R \cdot C \ =\  L_{R}\, Q(g,C) .
\end{eqnarray*}

We recall that if a semi-Riemannian manifold $(M,g)$, $n \geq 3$, 
is {\sl{locally symmetric}} then $\nabla R = 0$ on $M$. This  
implies the integrability condition ${\mathcal{R}}(X,Y ) \cdot R = 0$. In short, 
$R \cdot R = 0$. 
Semi-Riemannian manifolds satisfying this condition
are called 
{\sl semisymmetric} (see, e.g., {\cite[Chapter 8.5.3] {M7}} 
and {\cite[Chapter 1.6] {M6}}). 
Thus every semisymmetric manifold is pseudosymmetric.
The converse statement is not true (see, e.g., \cite{{DVV1991}}).
For examples of semisymmetric Roter spaces we refer to
{\cite[Example 5.4] {2015_DGHZ}}, {\cite[Section 3] {G5}}
and {\cite[Section 4] {Kow01}}.

A semi-Riemannian manifold $(M,g)$, $\dim M = n \geq 3$, 
is called {\sl Ricci-pseudosymmetric} 
if the tensors $R \cdot S$ and $Q(g,S)$ are linearly dependent 
at every point of $M$
(see, e.g., {\cite[Chapter 8.5.3] {M7}}, 
{\cite[Section 15.1] {CHEN-2021}}, \cite{DGHSaw}).
This is equivalent on ${\mathcal{U}}_{S}$ to 
$R \cdot S = L_{S}\, Q(g,S)$, 
where $L_{S}$ is some function on ${\mathcal U}_{S}$. 
Every warped product manifold $\overline{M} \times _{F} \widetilde{N}$
with a $1$-dimensional $(\overline{M}, \overline{g})$ manifold and
an $(n-1)$-dimensional Einstein 
semi-Riemannian manifold $(\widetilde{N}, \widetilde{g})$, $n \geq 3$, 
and a warping function $F$, is a Ricci-pseudosymmetric manifold
(see, e.g., {\cite[Section 1] {Ch-DDGP}}
and {\cite[Example 4.1] {2016_DGJZ}}).

It seems that (\ref{pseudo}) is the most important curvature condition 
between conditions called pseudosymmetry type conditions. 
We refer to 
\cite{{DGHP-TZ 1}, {DGHP-TZ 2}},
{\cite[Sections 3 and 4] {DGHSaw}},
{\cite[Section 1] {2016_DGJZ}},
{\cite[Section 2] {2018_DH}}
and
{\cite[Section 1] {2016_DP-TVZ}}
for a survey of results on semi-Riemannian manifolds 
satisfying such conditions. 
We mention that Roter spaces satisfy  (\ref{pseudo}) 
and some other pseudosymmetry type curvature conditions 
(see, e.g.,
{\cite[Section 4] {DGHP-TZ 1}}, {\cite[Theorem 6.7] {DGHSaw}}, 
{\cite[Theorem 3.2] {2016_DGJZ}}, {\cite[Theorem 2.1] {2018_DH}}).

\begin{lemma}
If $(M,g)$, $n \geq 4$, is a quasi-Einstein space satisfying
(\ref{quasi02}) on ${\mathcal U}_{S} \cap  {\mathcal U}_{C} \subset M$ 
then the following equation is satisfied on this set 
\begin{eqnarray}
g \wedge S^{2} + \frac{n-2}{2} \, S \wedge S - \kappa \, g \wedge S
+ \frac{\kappa ^{2} - \mathrm{tr}_{g} (S^{2})}{2(n-1)}
\, g \wedge g &=& 0 .
\label{2022.11.20.aa}
\end{eqnarray}
\end{lemma}
{\bf{Proof.}} Equation (\ref{quasi02}), 
in view of {\cite[Proposition 2.2] {G6}}, yields
\begin{eqnarray}
\frac{1}{2}\, S \wedge S &=& \alpha \, g \wedge S 
- \frac{\alpha ^{2}}{2}\, g \wedge g .
\label{2022.11.20.bb}
\end{eqnarray} 
From (\ref{2022.11.20.bb}), by suitable contractions, we get 
\begin{eqnarray}
S^{2} - \kappa \, S &=& - (n-2) \alpha \, S - \alpha \kappa \, g
+ (n-1) \alpha ^{2}\, g ,
\label{2022.11.20.cc}\\
\mathrm{tr}_{g} (S^{2}) - \kappa ^{2} &=&
- 2 (n-1)\alpha \kappa + n (n-1) \alpha ^{2} .
\label{2022.11.20.dd}
\end{eqnarray}
Now using (\ref{2022.11.20.bb}), (\ref{2022.11.20.cc}) and
(\ref{2022.11.20.dd}) we can easily check that 
(\ref{2022.11.20.aa}) 
holds on ${\mathcal U}_{S} \cap {\mathcal U}_{C}$.
Our lemma is thus proved.

\begin{lemma}
If $(M,g)$, $n \geq 4$, is a Roter space satisfying
(\ref{eq:h7a}) 
on ${\mathcal U}_{S} \cap  {\mathcal U}_{C} \subset M$ 
then the following equation is satisfied on this set 
\begin{eqnarray}
C &=& \frac{\phi}{n-2} 
\left(g \wedge S^{2} + \frac{n-2}{2} \, S \wedge S - \kappa \, g \wedge S
+ \frac{\kappa ^{2} - \mathrm{tr}_{g} (S^{2})}{2(n-1)}
\, g \wedge g \right) .
\label{2022.11.22.aa}
\end{eqnarray}
\end{lemma}
{\bf{Proof.}} Equation (\ref{eq:h7a}), 
by making use of (\ref{Weyl}), turns into
\begin{eqnarray}
C &=& \frac{\phi}{2}\, S\wedge S 
+ \left(\mu - \frac{1}{n-2} \right) g\wedge S + \frac{1}{2}
\left( \eta + \frac{\kappa }{(n-2) (n-1)} \right) g \wedge g .
\label{2022.11.22.bb}
\end{eqnarray}
From (\ref{2022.11.22.bb}), by suitable contractions, we get
\begin{eqnarray}
& &
\frac{\phi}{n-2}\, ( S^{2} - \kappa \, S)
- \left(\mu - \frac{1}{n-2} \right) S\nonumber\\
& & 
- \frac{\kappa }{n-2} \left(\mu - \frac{1}{n-2} \right) g
- \frac{n-1}{n-2} \left( \eta + \frac{\kappa }{(n-2) (n-1)} \right) g
\ =\ 0 ,
\label{2022.11.22.cc}\\
& &
\frac{\phi\, ( \kappa ^{2} - \mathrm{tr}_{g}(S^{2}))}{(n-2)(n-1)}
\ =\
- \frac{2 \kappa }{n-2} \left(\mu - \frac{1}{n-2} \right)
- \frac{n}{n-2} \left( \eta + \frac{\kappa }{(n-2) (n-1)} \right) .
\label{2022.11.22.dd}
\end{eqnarray}
Now (\ref{2022.11.22.bb}), (\ref{2022.11.22.cc}) 
and (\ref{2022.11.22.dd}) yield
\begin{eqnarray*}
C &=& \frac{\phi}{2}\, S\wedge S 
+ \left(\mu - \frac{1}{n-2} \right) g\wedge S + \frac{1}{2}
\left( \eta + \frac{\kappa }{(n-2) (n-1)} \right) g \wedge g\nonumber\\
& &
+ \frac{\phi}{n-2}\, g \wedge ( S^{2} - \kappa \, S)
- \left(\mu - \frac{1}{n-2} \right) g \wedge S \nonumber\\
& &
- \frac{\kappa}{n-2} \left(\mu - \frac{1}{n-2} \right) g \wedge g
- \frac{n-1}{n-2} \left( \eta + \frac{\kappa }{(n-2) (n-1)} \right) 
g \wedge g \nonumber\\
&=& \frac{\phi}{2}\, S\wedge S 
+ \frac{\phi}{n-2}\, g \wedge ( S^{2} - \kappa \, S)
- \frac{2 \kappa}{2 (n-2)} \left(\mu - \frac{1}{n-2} \right) g \wedge g
\nonumber\\
& &
- \frac{n}{2 (n-2)} \left( \eta + \frac{\kappa }{(n-2) (n-1)} \right) 
g \wedge g \nonumber\\
&=& \phi \left( \frac{1}{2}\, S\wedge S 
+ \frac{1}{n-2}\, g \wedge S^{2} - \frac{\kappa}{n-2}\, g \wedge S
+ \frac{  \kappa ^{2} - \mathrm{tr}_{g} (S^{2})}{2 (n-2)(n-1)}\, 
g \wedge g
\right) 
\end{eqnarray*}
and we get (\ref{2022.11.22.aa}), completing the proof.

Lemma 2.1 and Lemma 2.2 imply  
\begin{proposition}
If $(M,g)$, $n \geq 4$, is a semi-Riemannian manifold satisfying
(\ref{quasi02}) or (\ref{eq:h7a}) at every point of
${\mathcal U}_{S} \cap  {\mathcal U}_{C} \subset M$ 
then the following equation is satisfied on this set  
\begin{eqnarray}
\tau \, C &=&  
g \wedge S^{2} + \frac{n-2}{2} \, S \wedge S - \kappa \, g \wedge S
+ \frac{\kappa ^{2} - \mathrm{tr}_{g} (S^{2})}{2(n-1)}
\, g \wedge g ,
\label{2022.11.22.jj}
\end{eqnarray}
where $\tau$ is some function on ${\mathcal U}_{S} \cap {\mathcal U}_{C}$.
\end{proposition}

Proposition 2.3, {\cite[Theorem 7.1 (ii)] {2016_DGJZ}} and
{\cite[Theorem 4.1] {DK2003}} imply 
\begin{theorem} 
Let $\overline{M} \times _{F} \widetilde{N}$ 
be the warped product manifold 
with a $2$-dimensional semi-Riemannian manifold 
$(\overline{M},\overline{g})$,
an $(n-2)$-dimensional semi-Riemannian manifold 
$(\widetilde{N},\widetilde{g})$, $n \geq 4$, a warping function $F$, 
and let $(\widetilde{N},\widetilde{g})$ 
be a space of constant curvature when $n \geq 5$. 
Then (\ref{2022.11.22.jj}) holds on
${\mathcal U}_{S} \cap {\mathcal U}_{C} 
\subset \overline{M} \times _{F} \widetilde{N}$.
\end{theorem}

\begin{example}
(i) Let $\mathbb{S}^{k}(r^{-2})$ be a $k$-dimensional standard sphere      
of radius $r$ in $\mathbb{E}^{k+1}$, $k \geq 1$.
It is well-known that the Cartesian product 
$\mathbb{S}^{1}(r_{1}^{-2}) \times \mathbb{S}^{n-1}(r_{2}^{-2})$
of spheres $\mathbb{S}^{1}(r_{1}^{-2})$ and $\mathbb{S}^{n-1}(r_{2}^{-2})$, 
$n \geq 4$,
and more generally,
the warped product manifold 
$\mathbb{S}^{1}(r_{1}^{-2}) \times _{F} \mathbb{S}^{n-1}(r_{2}^{-2})$
of spheres $\mathbb{S}^{1}(r_{1}^{-2})$ and $\mathbb{S}^{n-1}(r_{2}^{-2})$,
$n \geq 4$, with a warping function $F$,
is a conformally flat manifold.
\newline
(ii)
As it was stated in {\cite[Example 3.2] {G5}},
the Cartesian product 
$\mathbb{S}^{p}(r_{1}^{-2}) \times \mathbb{S}^{n-p}(r_{2}^{-2})$
of spheres 
$\mathbb{S}^{p}(r_{1}^{-2})$ and $\mathbb{S}^{k}(r_{2}^{-2})$ 
such that 
$2 \leq p \leq n-2$ and
$(n-p-1) r_{1}^{2} \neq (p-1) r_{2}^{2}$
is a non-conformally flat and non-Einstein manifold
satisfying the Roter equation (\ref{eq:h7a}) on 
${\mathcal U}_{S} \cap  {\mathcal U}_{C} 
= \mathbb{S}^{p}(r_{1}^{-2}) \times \mathbb{S}^{n-p}(r_{2}^{-2})$. 
\newline
(iii) 
{\cite[Example 4.1] {DK2003}}
The warped product manifold 
$\mathbb{S}^{p}(r_{1}^{-2}) \times _{F} \mathbb{S}^{n-p}(r_{2}^{-2})$,
$2 \leq p \leq n-2$, with some special warping function $F$,
satisfies on ${\mathcal U}_{S} \cap  {\mathcal U}_{C} 
\subset \mathbb{S}^{p}(r_{1}^{-2}) \times \mathbb{S}^{n-p}(r_{2}^{-2})$ 
the Roter equation (\ref{eq:h7a}). Thus some 
warped product manifolds 
$\mathbb{S}^{2}(r_{1}^{-2}) \times _{F} \mathbb{S}^{n-2}(r_{2}^{-2})$ 
are Roter spaces.
\newline
(iv)  Properties of pseudosymmetry type of warped products 
with $2$-dimensional base manifold,
a warping function $F$,
and an $(n-2)$-dimensional fiber, 
$n \geq 4$,
assumed to be of constant curvature when $n \geq 5$, 
were determined in  
{\cite[Sections 6 and 7] {2016_DGJZ}}.
Evidently, 
warped product manifolds 
$\mathbb{S}^{2}(r_{1}^{-2}) \times _{F} \mathbb{S}^{n-2}(r_{2}^{-2})$,
$n \geq 4$, are such manifolds.
Let $g$, $R$, $S$, $\kappa$ and $C$ denote the metric tensor, 
the Riemann-Christoffel curvature tensor, the Ricci tensor, the scalar curvature 
and  
the Weyl conformal curvature tensor of 
$\mathbb{S}^{2}(r_{1}^{-2}) \times _{F} \mathbb{S}^{n-2}(r_{2}^{-2})$, respectively.
From {\cite[Theorem 7.1] {2016_DGJZ}}
it follows that on set $V$ of all points of 
${\mathcal U}_{S} \cap  {\mathcal U}_{C} 
\subset \mathbb{S}^{2}(r_{1}^{-2}) \times _{F} \mathbb{S}^{n-2}(r_{2}^{-2})$
at which the tensor $S^{2}$ is not a linear combination 
of the tensors $S$ and $g$,
the Weyl tensor $C$ is expressed by
\begin{eqnarray}
C &=& \lambda 
\left( g \wedge S^{2} + \frac{n-2}{2} \, S \wedge S - \kappa \, g \wedge S
+ \frac{\kappa ^{2} - \mathrm{tr}_{g} (S^{2}) }{2(n-1)} \, g \wedge g \right) ,
\label{2022.11.10.zz}
\end{eqnarray}
where $\lambda$ is some function on $V$.
This,  by (\ref{Weyl}), turns into
\begin{eqnarray*}
R &=& \lambda \,
g \wedge S^{2} + \frac{n-2}{2} \lambda \, S \wedge S 
+ \left( \frac{1}{n-2} - \kappa \lambda \right) g \wedge S\\
& &
+ \frac{1}{2(n-1)} \left( (\kappa ^{2} - \mathrm{tr}_{g} (S^{2}) ) \lambda
- \frac{\kappa }{n-2} \right) g \wedge g .
\end{eqnarray*}
Thus  (\ref{B001simply}) is satisfied on $V$.
Moreover,  (\ref{eq:h7a}) holds 
at all points of $({\mathcal U}_{S} \cap  {\mathcal U}_{C}) \setminus V$,
at which (\ref{quasi02}) is not satisfied.
From Lemma 2.2 it follows that (\ref{2022.11.10.zz})
holds at all points of ${\mathcal U}_{S} \cap  {\mathcal U}_{C} \subset
\mathbb{S}^{2}(r_{1}^{-2}) \times _{F} \mathbb{S}^{n-2}(r_{2}^{-2})$, $n \geq 4$,
at which (\ref{quasi02}) is not satisfied. Finally, in view of Theorem 2.4, 
we can state that (\ref{2022.11.22.jj}) holds on
${\mathcal U}_{S} \cap  {\mathcal U}_{C}$.
\end{example}

\begin{remark}
We define on a semi-Riemannian manifold $(M,g)$, $\dim M = n \geq 3$, 
the $(0,4)$-tensor $E$ by {\cite[eq. (1.2)] {DGHP-TZ 1}} 
\begin{eqnarray*}
E &=& 
g \wedge S^{2} + \frac{n-2}{2} \, S \wedge S - \kappa \, g \wedge S
+ \frac{\kappa ^{2} - \mathrm{tr}_{g} (S^{2}) }{2(n-1)} \, g \wedge g .
\end{eqnarray*}
Evidently, we can present now 
(\ref{2022.11.20.aa}), (\ref{2022.11.22.aa}) and (\ref{2022.11.10.zz}) 
as follows
$E = 0$, $C = (\phi / (n-2)) E$ and $C = \lambda E$, respectively.
We refer to \cite{DGHP-TZ 1} for further results on semi-Riemannian manifolds
and in particular, hypersurfaces in spaces of constant curvature, satisfying
curvature conditions involving the tensor $E$.
\end{remark}

\section{Hypersufaces satisfying ${\mathcal A}^{2} = \alpha {\mathcal A} + \beta I\!d$.}

Let $M$ be a connected hypersurface isometrically immersed in
a semi-Riemannian space of constant curvature
$N_{s}^{n+1}(c)$, $n \geq 3$, 
with signature $(s, n+1-s)$, 
where 
$c = \frac{\widetilde{\kappa}}{n(n+1)}$ and
$\widetilde{\kappa}$ is the scalar curvature of the ambient space.
Let $g$ be the metric tensor induced on $M$ from the metric tensor 
of $N_{s}^{n+1}(c)$.
We denote by  
$R$, $S$, $\kappa$ and $C$ 
the Riemann-Christoffel curvature tensor,
the Ricci tensor, the scalar curvature and the Weyl conformal curvature tensor
of $M$, respectively. 
Further, let
$H$ and ${\mathcal A}$ be the second fundamental tensor and the shape operator
of $M$,
respectively. We have $H(X,Y) = g( {\mathcal A}X,Y)$,
for any vector fields $X,Y$ tangent to $M$.
The $(0,2)$-tensors $H^{2}$, $H^{3}$ and $H^{4}$ are defined by
\begin{eqnarray*}
H^{2}(X,Y) \ =\ H( {\mathcal A}X,Y) ,\ \ \
H^{3}(X,Y) & =& H^{2}( {\mathcal A}X,Y) ,\ \ \
H^{4}(X,Y) \ =\ H^{3}( {\mathcal A}X,Y) ,
\end{eqnarray*}
respectively, 
where $X,Y$ are vectors fields tangent to $M$.

From (\ref{realC5}), by contraction with $g^{ij}$, we obtain
\begin{eqnarray}
S_{hk} - (n-1) c\, g_{hk} 
&=&  \varepsilon \, ( \mathrm{tr} (H)\, H_{hk} - H^{2}_{hk}) ,
\label{realC6}\\
\kappa - n (n-1) c  
&=&  \varepsilon \, ( ( \mathrm{tr} (H))^{2} - \mathrm{tr} (H^{2}) ) ,
\label{realC6bb}
\end{eqnarray}
where 
$S_{hk}$ and $g^{hk}$ are the local components of the Ricci tensor $S$ 
and the tensor $g^{-1}$ of $M$, res\-pectively,  
and
$\mathrm{tr} (H) = g^{ij}H_{ij}$, 
$H^{2}_{hk} = g^{ij} H_{hi}H_{jk}$,
$\mathrm{tr} (H^{2}) = g^{ij}H^{2}_{ij}$. 
Further, from (\ref{realC6}) we get immediately
\begin{eqnarray}
\ \ \ \
S^{2}_{hk} - 2 (n-1) c\, S_{hk} + (n-1)^{2} c^{2}\, g_{hk} 
\ =\ 
H^{4}_{hk} - 2 \mathrm{tr} (H)\, H^{3}_{hk} 
+ 
(\mathrm{tr} (H))^{2} \, H^{2}_{hk} ,
\label{realC7}
\end{eqnarray}
where $S^{2}_{hk} = g^{ij} S_{hi}S_{jk}$,
$H^{3}_{hk} = g^{ij} H^{2}_{hi}H_{jk}$ and  
$H^{4}_{hk} = g^{ij} H^{3}_{hi}H_{jk}$,
are the local components of the $(0,2)$-tensors
$S^{2}$, $H^{3}$ and $H^{4}$, respectively.

Let $M$ be a non-Einstein hypersurface in  
a semi-Riemannian space of constant curvature $N_{s}^{n+1}(c)$, $n \geq 3$.
We assume that at a point $x \in {\mathcal{U}}_{S} \subset M$ 
the tensor $H^{2}$ is a linear combination of the tensors 
$g$ and $H$, i.e., 
$x \in {\mathcal{U}}_{S} \setminus {\mathcal U}_{H}$.
Evidently, the tensor $S - (\kappa /n)\, g$ is a non-zero tensor at $x$
and (\ref{2021.10.25.a2}) holds at $x$. 
Now (\ref{realC6}), by (\ref{2021.10.25.a2}), turns into
\begin{eqnarray} 
S &=& \varepsilon ( \mathrm{tr} (H) - \alpha )\, H 
+ ( (n-1)c - \varepsilon \beta )\, g .
\label{2021.10.25.a3}
\end{eqnarray} 
This together with (\ref{2021.10.25.a2}) yields 
\begin{eqnarray} 
S^{2} &=& ( \mathrm{tr} (H) - \alpha )^{2} \, H^{2} 
+ 2 \varepsilon ( \mathrm{tr} (H) - \alpha ) 
( (n-1)c - \varepsilon \beta )\, H 
+ ( (n-1)c - \varepsilon \beta )^{2}\, g \nonumber\\
&=&
( \mathrm{tr}\, (H) - \alpha )^{2} \, 
(  \alpha \, H + \beta \, g)\nonumber\\
& &
+ 2 \varepsilon ( \mathrm{tr} (H) - \alpha ) 
( (n-1)c - \varepsilon \beta ) \, H
+ 
( (n-1)c - \varepsilon \beta )^{2} \, g 
\nonumber\\
&=& ( \mathrm{tr} (H) - \alpha )
( \alpha \, ( \mathrm{tr}\, (H) - \alpha ) 
+ 2 \varepsilon  
( (n-1)c - \varepsilon \beta ) ) \, H \nonumber\\
& &
+
( \beta \, ( \mathrm{tr}\, (H) - \alpha )^{2} 
+
( (n-1)c - \varepsilon \beta )^{2} )\, g .
\label{2021.10.25.a4}
\end{eqnarray} 
If we would have
$\mathrm{tr} (H) - \alpha = 0$ at $x$
then
(\ref{2021.10.25.a3}) leads to 
$S = ( (n-1)c - \varepsilon \beta )\, g$, 
and as a consequence we get
$S = ( \kappa / n)\, g$, a contradiction.
Therefore
\begin{eqnarray}
\mathrm{tr} (H) - \alpha \neq 0  
\label{2021.10.25.a55}
\end{eqnarray}
at this point.
Now from
(\ref{2021.10.25.a3}), by (\ref{2021.10.25.a55}), we obtain
\begin{eqnarray} 
H &=&  \varepsilon ( \mathrm{tr} (H) - \alpha )^{-1}
( S - ( (n-1)c - \varepsilon \beta )\, g ).
\label{2021.10.25.a5}
\end{eqnarray} 
Next, using 
(\ref{2021.10.25.a4}), (\ref{2021.10.25.a55})
and
(\ref{2021.10.25.a5})
we can check that $M$ is partially Einstein at $x$.
We also note that (\ref{realC5}), by making use of 
(\ref{2021.10.25.a5}), turns into
\begin{eqnarray} 
\ \ \ \ \ \
R &=& \frac{ \varepsilon }{2 ( \mathrm{tr} (H) - \alpha )^{2} }\,
( S - ( (n-1)c - \varepsilon \beta )\, g ) \wedge
( S - ( (n-1)c - \varepsilon \beta )\, g ) + \frac{c}{2}\, g \wedge g .
\label{2021.10.25.a6}
\end{eqnarray} 
Thus we have
\begin{proposition}
Let $M$ be a hypersurface in a semi-Riemannian space of constant curvature
$N_{s}^{n+1}(c)$, $n \geq 3$. Then
(\ref{2021.10.25.a4}), (\ref{2021.10.25.a55}) and (\ref{2021.10.25.a5})
are satisfied  
at every point 
of ${\mathcal{U}}_{S} \setminus {\mathcal U}_{H} \subset M$.
\end{proposition}
\begin{proposition}
$\mathrm{ ( }$cf. {\cite[Theorem 3.3] {2018_DGZ}}, 
{\cite[Proposition 3.3] {G5}}$\mathrm{ ) }$
Let $M$ be a hypersurface in a semi-Riemannian space of constant curvature
$N_{s}^{n+1}(c)$, $n \geq 4$. 
Then (\ref{2021.10.25.a6}) is satisfied
at every point of  
$({\mathcal{U}}_{S} \cap {\mathcal{U}}_{C}) 
\setminus {\mathcal U}_{H}   \subset M$, i.e.,
an equation of the form (\ref{eq:h7a}).
\end{proposition}

Let $M$ be a hypersurface in a Riemannian space of constant curvature $N^{n+1}(c)$,
$n \geq 3$. We assume that $M$ has at every point exactly two distinct 
principal curvatures $\lambda_{1}$ and $\lambda_{2}$, with multiplicities 
$p$ and $n-p$, respectively, where $1 \leq p \leq n-1$. 
Thus we have on $M$: (\ref{2021.10.25.a2}) and
\begin{eqnarray} 
\ \ \ 
\alpha \ =\ \lambda_{1} + \lambda_{2} ,\ \ \
\beta \ =\ - \lambda_{1} \lambda_{2},\ \ \
\mathrm{tr} (H) - \alpha \ =\ (p-1) \lambda_{1} + (n-p-1) \lambda_{2} . 
\label{2022.09.16.bb}  
\end{eqnarray}
Further, from (\ref{realC6}) it follows that the eigenvalues 
$\kappa_{1}$ and $\kappa_{2}$ of the Ricci operator ${\mathcal{S}}$ 
of $M$ are expressed by 
\begin{eqnarray}
\kappa_{1} \ =\ - \lambda_{1}^{2} 
+ \mathrm{tr}(H) \lambda_{1} + (n-1) c ,\ \ \
\kappa_{2} \ =\ - \lambda_{2}^{2} 
+ \mathrm{tr}(H) \lambda_{2} + (n-1) c .
\label{2022.01.12.a2}
\end{eqnarray} 
Now from 
(\ref{2022.09.16.bb}) and (\ref{2022.01.12.a2})
we get immediately
\begin{eqnarray}
\kappa_{1} - \kappa_{2} &=&
- \lambda_{1}^{2} + \lambda_{2}^{2} 
+ \mathrm{tr}(H) (\lambda_{1} - \lambda_{2})
\ =\ 
(\lambda_{1} - \lambda_{2})
(\mathrm{tr}(H) - \alpha )\nonumber\\
&=& 
(\lambda_{1} - \lambda_{2})
( (p-1) \lambda_{1} + (n-p-1) \lambda_{2} ) .
\label{2022.01.12.a3}
\end{eqnarray} 
Evidently, if at every point of $M$ we have
$\mathrm{tr}(H) - \alpha \neq 0$,
or, equivalently,
(\ref{2022.09.23.cc}) is satisfied,
then 
${\mathcal{U}}_{S} = M$
and
$S^{2} = ( \kappa_{1} + \kappa_{2})\, S - \kappa_{1} \kappa_{2} \, g$
on $M$, 
i.e., such a hypersurface is partially Einstein.
In addition, we assume that
$n \geq 4$ and $2 \leq p \leq n-2$.
Thus $M$ is a non-quasi-umbilical at every point, 
and in a consequence,
$M$ is a non-conformally flat hypersurface. Thus we have
${\mathcal{U}}_{S} \cap {\mathcal{U}}_{C} = M$. 
Now from Proposition 3.2 it follows that
$M$ is a Roter hypersurface. Precisely, we have
\begin{theorem}
Let $M$ be a hypersurface in a Riemannian space of constant curvature
$N^{n+1}(c)$, $n \geq 4$, such that at every point $M$ has exactly two distinct 
principal curvatures, $\lambda_{1}$ with multiplicity $p$ and
$\lambda_{2}$ with multiplicity $n-p$, $2 \leq p \leq n-2$.
If (\ref{2022.09.23.cc}) is satisfied 
at every point of $M$
then $M = {\mathcal{U}}_{S}  \cap {\mathcal{U}}_{C}$
and an equation of the form (\ref{eq:h7a}) holds on $M$,
i.e., $M$ is a Roter space. In particular, 
$M$ is a non-quasi-Einstein partially Einstein space.
\end{theorem}
\begin{example} {\cite[Example 1.1] {Wu-2011}}
We present now a family 
of non-Einstein and non-conformally flat hypersurfaces $M$ in
$N^{n+1}(c)$, $n \geq 4$,
having at every point exactly two distinct principal curvatures.
Let 
\begin{eqnarray*}
N^{n+1}(c) \ =\ \mathbb{S}^{n+1}(c) &=& 
\{ y \in \mathbb{R}^{n+2}\, : \, < y,y > \ =\ \frac{1}{c} ,\, c > 0 \},
\end{eqnarray*}
where $< \cdot , \cdot >$ is the standard inner product on 
$\mathbb{R}^{n+2}$, $n \geq 3$. 
Let $M_{p,n-p} (c,t)$ be hypersurfaces in $\mathbb{S}^{n+1}(c)$ 
defined by
\begin{eqnarray}
M_{p,n-p} (c,t)  &=&
\mathbb{S}^{p}( c / \sin ^{2} t ) 
\times 
\mathbb{S}^{n-p} ( c / \cos ^{2} t ) ,
\end{eqnarray}
where $t \in ( 0, \frac{\pi }{2})$ and $1 \leq p \leq n-1$.
The hypersurface $M_{p,n-p} (c,t)$ is called  the {\sl{Clifford hypersurface}} 
(see, e.g., {\cite[Example 1.1] {Wu-2011}}). 
The Clifford hypersurface is also called the {\sl{generalized Clifford torus}}
(see, e.g., {\cite[p. 13] {BNO-2023}}).
The principal curvatures
$\lambda_{1}$ and $\lambda_{2}$ of
every hypersurface $M_{p,n-p} (c,t)$ are the following:
$\lambda_{1}\, =\, \sqrt{c}\, \cot t > 0$  
and
$\lambda_{2}\, =\, - \sqrt{c}\, \tan t < 0$.
Hypersurfaces 
$M_{p,n-p} (c,t)$ are non-conformally flat,
provided that  $n \geq 4$ and $2 \leq p \leq n-2$. In addition,
from (\ref{2022.01.12.a3}) it follows that a hypersurface
$M_{p,n-p} (c,t)$ is a non-Einstein hypersurface 
if and only if 
(\ref{2022.09.23.cc}) is satisfied 
at every point of $M_{p,n-p} (c,t)$.
The last condition is satisfied when
$t \neq \arctan \sqrt{ ( p-1 ) / ( n -p -1) }$. 
\end{example}
\begin{remark}
Hypersurfaces in Riemannian spaces of constant curvature
having at every point two distinct principal curvatures 
satisfy some curvature conditions of pseudosymmetry type, 
see, e.g., {\cite[Section 5] {2016_DGHZhyper}}
and {\cite[Section 3] {2018_DGZ}}. 
\end{remark}


\section{Hypersufaces satisfying ${\mathcal A}^{3} = \phi {\mathcal A}^{2} + \psi {\mathcal A} + \rho I\!d$. The general case.}

Let $(M,g)$, $\dim M \geq 4$, be a hypersurface in an $(n+1)$-dimensional
semi-Riemannian manifold $(N,\tilde{g})$, 
where the metric tensor $g$ of $M$ is induced 
by the metric tensor $\tilde{g}$ of the ambient space $N$.
The hypersurface $(M,g)$ is said to be {\sl{quasi-umbilical}}
(see, e.g., \cite{{DV_1991_GTS}, {G6}})
resp., $2$-{\sl{quasi-umbilical}}
(see, e.g., \cite{{DeVerYap}, {G6}}) 
at a point $x \in M$
if $\mathrm{rank} (H - \alpha \, g) = 1$, 
resp., $\mathrm{rank} (H - \alpha \, g) = 2$ holds at $x$,
for some $\alpha \in \mathbb{R}$.
If at a point  $x \in M$ we have $H =  \alpha \, g$,  
for some $\alpha \in \mathbb{R}$, then $M$ is called 
{\sl{umbilical}} at $x$.

If $(M,g)$, $\dim M \geq 4$, is a hypersurface in an $(n+1)$-dimensional
Riemannian manifold $(N,\tilde{g})$ then the above presented definitions of 
a  quasi-umbilical and $2$-quasi-umbilical points of the hypersurface 
$M$ are equivalent to the following definitions.
Namely,
according to 
\cite{{CHEN-2000}, {CHEN-VERSTRAELEN-1976}, {OU-CHEN-2020}, {PRV_2002},
{V2}} (see also \cite{DGP-TV02}),
the hypersurface $M$ is said to be {\sl{quasi-umbilical}},
resp., $2$-{\sl{quasi-umbilical}}, at a point $x \in M$
if it has at $x$ a principal curvature with multiplicity $n-1$, 
resp., $n-2$,
i.e., when the principal curvatures of $M$ at $x$ are the following:
$\lambda _{1} \neq \lambda _{2} = \lambda _{3} = \ldots = \lambda _{n}$, 
resp.,
$\lambda _{1}$, $\lambda _{2}$, 
$\lambda _{3} = \lambda _{4} = \ldots = \lambda _{n}$,
with  
$\lambda _{1} \neq \lambda _{3}$ and $\lambda _{2} \neq \lambda _{3}$.

We recall that a hypersurface $(M,g)$, $\dim M = n \geq 4$,
in a semi-Riemannian conformally flat manifold $(N,\tilde{g})$,
where $g$ is induced by $\tilde{g}$,
is quasi-umbilical
at a point $x \in M$ if and only if its Weyl conformal curvature tensor 
$C$ vanishes at this point {\cite[Theorem 4.1] {DV_1991_GTS}}.
Thus $M$ is non-quasi-umbilical at a point $x \in M$
if and only if its Weyl conformal curvature tensor $C$ 
is non-zero at $x$, i.e., $x \in \mathcal{U}_{C} \subset M$.

Let now $M$ be a hypersurface in $N_{s}^{n+1}(c)$, $n \geq 3$,
with the tensor $H$ satisfying (\ref{real702})  on ${\mathcal U}_{H}$.
From (\ref{real702}) we get on ${\mathcal U}_{H}$
\begin{eqnarray}
H^{4} &=& \phi \, H^{3} + \psi \, H^{2} + \rho \, H 
\ =\ ( \phi ^{2} + \psi ) \, H^{2} 
+ (\phi \psi + \rho)\, H + \phi \rho \, g .
\label{real703} 
\end{eqnarray}
Now (\ref{realC7}), with the aid of (\ref{real702}),  (\ref{realC6}) 
and (\ref{real703}), turns into 
\begin{eqnarray}
A &=& \tau \, H ,
\label{real70456} 
\end{eqnarray} 
where the $(0,2)$-tensor $A$ and the function $\tau$ are defined by
\begin{eqnarray}
A &=&  S^{2} - \rho _{1}\, S - \rho _{0}\, g ,
\label{real707}\\
\tau &=& \rho 
+ ( \phi - \mathrm{tr} (H) ) 
( \psi  + \mathrm{tr} (H) ( \phi - \mathrm{tr} (H) ) )  ,
\label{real705}
\end{eqnarray}
respectively, where
\begin{eqnarray}
\rho _{0} &=& 
- ( (n-1)^{2} c^{2} 
- \varepsilon  ( \psi + ( \phi - \mathrm{tr} (H) )^{2} )  (n-1) c
+ \rho (\phi - 2 \mathrm{tr} (H) ) ) ,
\label{real707rho0}\\  
\rho _{1} &=& 2 (n-1) c 
- \varepsilon ( \psi +  ( \phi - \mathrm{tr} (H) )^{2} )  .
\label{real707rho1}
\end{eqnarray}
We note that (\ref{realC5}), 
by (\ref{real70456}) and (\ref{real707}),  
leads to
\begin{eqnarray}
\tau ^{2}\, R &=& \frac{\varepsilon}{2}\, A \wedge A 
+ \frac{\tau ^{2} c}{2}\, g \wedge g
\label{2022.01.11.a1}
\end{eqnarray}
and (\ref{2022.09.19.aa}).
From 
(\ref{realC6}), 
(\ref{real70456}) and (\ref{real707}) it follows that
\begin{eqnarray}
\rho _{0} &=& \frac{1}{n}
\left(
\mathrm{tr} (S^{2}) - \rho_{1} \kappa - \tau \mathrm{tr}(H) \right) ,
\nonumber\\
S^{2} - \rho _{1}\, S - \tau \, H
&=& 
\frac{1}{n}
\left( \mathrm{tr} (S^{2}) - \rho _{1}\, \kappa 
- \tau \, \mathrm{tr} (H) \right) g ,
\label{2021.12.17.a6}
\end{eqnarray}
\begin{eqnarray}
Q(g,S^{2})  &=& \rho _{1}\, Q(g,S) + \tau \, Q(g,H) ,
\label{2021.12.17.a5}\\
\tau \, Q(H,S)  &=& 
Q(\tau \, H, S) \ =\ Q(A,S) 
\ =\  Q( S^{2} - \rho _{1}\, S - \rho _{0}\, g  ,S) 
\nonumber\\
&=&
Q(S^{2}, S) - \rho _{0} \, Q(g ,S) 
\ = \ -  Q(S, S^{2}) - \rho _{0} \, Q(g ,S) ,
\label{2022.01.07.a7}\\
\tau \, Q(g,H)  &=& 
Q(g,\tau \, H) \ =\ Q(g, A) 
\ =\
Q(g, S^{2} - \rho _{1}\, S - \rho _{0}\, g )
\nonumber\\
&=&
Q(g, S^{2} ) - \rho _{1}\, Q(g,S) ,
\label{2022.01.07.77}
\end{eqnarray}
\begin{eqnarray}
\tau \, Q(H,H^{2})  &=& 
\tau \, Q(H, \mathrm{tr} (H)\, H -\varepsilon \, S 
+ (n-1) \varepsilon c\, g ) \nonumber\\
&=& \varepsilon \tau \, Q(H, - S + (n-1) c\, g )
\ =\
\varepsilon  \, Q( A , - S + (n-1) c\, g )
\nonumber\\
&=& \varepsilon 
\left( Q(S, S^{2})
+ (\rho_{0} + (n-1)c \rho_{1})\, Q(g, S)
- (n-1) c\, Q(g, S^{2})
\right) .
\label{2022.01.07.88}
\end{eqnarray}
Moreover, from (\ref{real70456}) and (\ref{real707}) we also obtain
\begin{eqnarray}
S^{2}_{ij} &=& \rho _{1}\, S_{ij} + \rho _{0}\, g_{ij} + \tau\, H_{ij} ,
\label{new01-real706} \\
S^{3}_{ij} &=& \rho _{1}\, S^{2}_{ij} + \rho _{0}\, S_{ij} 
+ \tau H_{ir}g^{rs}S_{sj} , 
\label{new02-real706} 
\end{eqnarray}
where $S^{3}_{ij} = S_{ih}g^{hk}S^{2}_{kj}$ are the local components 
of the $(0,2)$-tensor $S^{3}$. 
Now applying 
(\ref{real702}), (\ref{realC6}) and (\ref{real70456}) into (\ref{new02-real706}) 
we get 
\begin{eqnarray*}
S^{3}_{ij} &=& \rho _{1}\, S^{2}_{ij} + \rho _{0}\, S_{ij} 
+ \tau H_{ir} g^{rs} ( 
\varepsilon \, (\mathrm{tr} (H)\, H_{sj} - H^{2}_{sk} ) 
+ (n-1) c\, g_{sj} ) \nonumber\\
&=& \rho _{1}\, S^{2}_{ij} + \rho _{0}\, S_{ij} 
+ \tau (
\varepsilon \, ( \mathrm{tr} (H)\, H^{2}_{ij} - H^{3}_{ij}) 
+ (n-1) c\, H_{ij} ) \nonumber\\
&=& \rho _{1}\, S^{2}_{ij} + \rho _{0}\, S_{ij} 
+ \varepsilon \, \tau ( \mathrm{tr} (H)\, H^{2}_{ij} - \phi H^{2}_{ij} 
- \psi \, H_{ij} - \rho \, g_{ij} 
+ \varepsilon (n-1) c\, H_{ij} ) \nonumber\\
&=& \rho _{1}\, S^{2}_{ij} + \rho _{0}\, S_{ij} 
+ \varepsilon \, \tau ( ( \mathrm{tr} (H) - \phi ) H^{2}_{ij} 
+ (- \psi + \varepsilon  (n-1) c)\, H_{ij} - \rho\, g_{ij} ) ,
\end{eqnarray*}
i.e.,
\begin{eqnarray}
S^{3} &=& \rho _{1}\, S^{2} + \rho _{0}\, S 
+ \varepsilon  \tau ( ( \mathrm{tr} (H) - \phi )\, H^{2} 
+ ( \varepsilon (n-1) c - \psi )\, H - \rho \, g) .
\label{new04-real7077} 
\end{eqnarray}

Let $M$ be a hypersurface in $N_{s}^{n+1}(c)$, $n \geq 4$,
with the tensor $H$ satisfying (\ref{real702})  
on ${\mathcal U}_{H} \subset M$.
We define on ${\mathcal U}_{H}$ the following functions
({\cite[eq. (34)] {Saw-2004}}, {\cite[Section 4] {Saw-2007}}):
\begin{eqnarray}
\beta_{1} &=& \varepsilon (\phi - \mathrm{tr}(H)) ,\nonumber\\
\beta_{2} &=& - \frac{\varepsilon }{n-2} 
( \phi \, (2 \mathrm{tr}(H) - \phi)
- (\mathrm{tr}(H))^{2} - \psi - (n-2) \varepsilon \mu ) \nonumber\\
&=& 
\mu + 
\frac{\varepsilon }{n-2} 
\left( \beta_{1}^{2} + \psi \right) ,\nonumber\\
\beta_{3} &=& \varepsilon \mu \mathrm{tr}(H)
+ \frac{1}{n-2}\, (\psi \, (2 \mathrm{tr}(H) - \phi)
 + (n-3) \rho ) ,
\label{ZZ23}\\
\beta_{4} &=& \beta_{3} -  \varepsilon \beta_{2} \mathrm{tr}(H)
+ \frac{ (n-1) \widetilde{\kappa} \beta_{1} }{n(n+1)}  ,\nonumber\\
\beta_{5} &=&
\frac{\kappa }{n-1}
+ \varepsilon \psi 
- \frac{ (n^{2} - 3n + 3) \widetilde{\kappa} }{n(n+1)}
+ \beta_{1} \mathrm{tr}(H) ,
\nonumber\\
\beta_{6} &=& \beta_{2} 
- \frac{ (n-3) \widetilde{\kappa} }{n(n+1)} ,\nonumber
\end{eqnarray}
\begin{eqnarray}  
\mu
&=&
\frac{1}{n-2}
\left( \frac{\kappa }{n-1}
- \frac{\widetilde{\kappa} }{n+1}
\right) ,
\label{newZZ01}
\end{eqnarray}
where $\phi$, $\psi $ and $\rho $ 
are defined by (\ref{real702}).
We also have on ${\mathcal U}_{H}$
{\cite[eq. (24)] {Saw-2007}}
\begin{eqnarray}
R \cdot S
&=&
\frac{\widetilde{\kappa}}{n (n+1)}\, Q(g,S)
+ \rho \, Q(g,H)
- \varepsilon \beta _{1}\, Q(H, H^{2}) 
\label{900r1}
\end{eqnarray}
and
{\cite[eqs. (52), (45), (46), (47)] {Saw-2004}}
\begin{eqnarray}
C \cdot S
&=&
\beta _{1}\, Q(H,S)
+ \beta _{2}\, Q(g,S)
+ \beta _{4} \, Q(g,H) ,
\label{newZZ04}
\end{eqnarray}
\begin{eqnarray}
(n-2)\, R \cdot C  &=& (n-2)\, Q(S,R) 
- \frac{(n-2)^{2} \widetilde{\kappa }}{n (n+1)}\, Q(g,R)
- \frac{(n-3) \widetilde{\kappa} }{2 n (n+1)}\, Q(S,g \wedge g)
\nonumber\\
& &
+ \frac{\rho }{2} \, Q(H, g \wedge g)
+ \varepsilon \beta _{1} \, g \wedge Q(H,H^{2}) ,
\label{ZZ1}
\end{eqnarray}
\begin{eqnarray}
(n-2)\, C \cdot R &=&
(n-3)\, Q(S,R) + \varepsilon \beta _{1} \, H \wedge Q(g,H^{2})
- \frac{(n - 3) \widetilde{\kappa }}{2 n (n+1)}\, Q(S, g \wedge g )
\nonumber\\
& &
+
 \left( \frac{\kappa }{n-1} +  \varepsilon  \psi
- \frac{(n^{2} - 3n + 3) \widetilde{\kappa }}{n(n+1)}
\right) Q(g,R) ,
\label{ZZ2}
\end{eqnarray}
\begin{eqnarray}
(n-2)\, C \cdot C &=&
(n-3)\, Q(S,R) 
+ \beta _{5}\, Q(g,R) + \frac{\beta _{6}}{2} Q(S, g \wedge g )
\nonumber\\
& &
+ \beta _{1} \, Q(S, g \wedge H) 
+ \frac{\beta _{4}}{2} \, Q(H, g \wedge g) , 
\label{2021.11.30.a5}
\end{eqnarray} 
where $\beta_{1}, \ldots , \beta_{6}$
are defined by (\ref{ZZ23}).

From (\ref{900r1}), by making use of 
(\ref{2022.01.07.77}) and (\ref{2022.01.07.88}), we get
\begin{eqnarray}
\tau \, R \cdot S
&=&
- \beta _{1} \, Q(S, S^{2}) 
+ ( \rho + (n-1) c \beta_{1} )\, Q(g, S^{2})\nonumber\\
& &
+ ( \tau c -  \rho  \rho _{1} 
- \beta _{1} ( \rho _{0} + (n-1) c  \rho _{1}) )\, Q(g,S) .
\label{2022.01.10.99}
\end{eqnarray} 
We also note that (\ref{newZZ04}), 
by (\ref{2021.12.17.a5}) and (\ref{2022.01.07.a7}), turns into
\begin{eqnarray}
\tau \, C \cdot S
&=&
- \beta _{1}\, Q(S,S^{2})
+ \beta _{4} \, Q(g,S^{2})
+
(\beta _{2} \tau - \beta _{1} \rho _{0} - \beta _{4} \rho _{1} )
\, Q(g,S)  .
\label{2022.01.10.j}
\end{eqnarray} 
Thus we have
\begin{proposition}
Let $M$ be a hypersurface in a semi-Riemannian space of constant curvature
$N_{s}^{n+1}(c)$, $n \geq 4$,
with the tensor $H$ satisfying (\ref{real702}) 
on ${\mathcal U}_{H} \subset M$.
Then 
(\ref{2022.09.19.aa}), 
(\ref{real70456}), 
(\ref{real705}) and
(\ref{2022.01.11.a1})--(\ref{2022.01.10.j})
hold
on ${\mathcal U}_{H}$.
\end{proposition}

\begin{remark} Conditions (\ref{real70456}), 
(\ref{real705}) and
(\ref{2022.01.11.a1})
were already presented in the last section of a preliminary version of 
\cite{2021-DGH}. Precisely, see eqs. (8.6)--(8.9) 
in Section 8 of arXiv:1911.02482v2 [math.DG] 27 Jan 2020.
However, that section was not included to \cite{2021-DGH}.
\end{remark}


\section{Hypersufaces satisfying ${\mathcal A}^{3} = \phi {\mathcal A}^{2} + \psi {\mathcal A} + \rho I\!d$. The partially Einstein case.}

Let $M$ be a hypersurface in $N_{s}^{n+1}(c)$, $n \geq 3$,
with the tensor $H$ satisfying (\ref{real702}) 
on ${\mathcal U}_{H} \subset M$.
We assume that the hypersurface $M$ is a partially Einstein 
at a point
$x \in {\mathcal U}_{H}$, i.e., let
$S^{2} = \widetilde{\alpha }\, S + \widetilde{\beta }\, g$ at $x$, 
where $\widetilde{\alpha }, \, \widetilde{\beta} \in \mathbb{R}$.
Now (\ref{real70456}) turns into 
\begin{eqnarray}
(\widetilde{\alpha } - \rho _{1}) \, S 
+ (\widetilde{\beta } - \rho _{0} ) \, g  &=& \tau \, H  .
\label{2022.01.15.c} 
\end{eqnarray}
If we would have $\widetilde{\alpha } \neq \rho _{1}$ then
$S = (\widetilde{\alpha } - \rho _{1})^{-1} 
\, (\tau \, H - ( \widetilde{\beta } - \rho _{0} ) \, g)$
at $x$. This by (\ref{realC6}) leads to an equation of the form
(\ref{2021.10.25.a2}), a contradiction. 
Therefore $\widetilde{\alpha } = \rho _{1}$ at $x$.
Now (\ref{2022.01.15.c}) reduces to
$(\widetilde{\beta } - \rho _{0} ) \, g  = \tau \, H$.
From this it follows that  
$\widetilde{\beta } = \rho _{0}$ and $\tau = 0$ at $x$.
Thus we have
\begin{theorem}
Let $M$ be a hypersurface in a semi-Riemannian space of constant curvature
$N_{s}^{n+1}(c)$, $n \geq 3$,
with the tensor $H$ satisfying (\ref{real702})
on ${\mathcal U}_{H} \subset M$.
If $M$ is partially Einstein at every point of ${\mathcal U}_{H}$
then 
$S^{2} = \rho _{1}\, S + \rho _{0}\, g$ and $\tau = 0$, i.e.,
\begin{eqnarray*}
\rho 
+ ( \phi - \mathrm{tr} (H) ) 
( \psi  + \mathrm{tr} (H) ( \phi - \mathrm{tr} (H) ) ) &=& 0 
\end{eqnarray*}
on ${\mathcal U}_{H}$, where $\rho _{0}$, $\rho _{1}$ and $\tau $ are defined by
(\ref{real707rho0}), (\ref{real707rho1}) 
and (\ref{real705}), respectively.
\end{theorem}


Let $M$ be a hypersurface in a Riemannian space of constant curvature
$N^{n+1}(c)$, $n \geq 3$,
having at every point $x \in {\mathcal U}_{H} \subset M$
three distinct principal curvatures 
$\lambda _{0}$, $\lambda _{1}$ and $\lambda _{2}$,
with multiplicities $n_{0}$, $n_{1}$ and $n_{2}$, respectively,
where $n_{0} + n_{1} + n_{2} = n$.
Evidently,  
(\ref{real702}) holds on ${\mathcal U}_{H}$ with 
\begin{eqnarray}
\phi \ =\ \lambda _{0} + \lambda _{1} + \lambda _{2},\ \ \ 
\psi \ =\  - \lambda _{0}
 ( \lambda _{1} + \lambda _{2}) - \lambda _{1} \lambda _{2},\ \ \
\rho \ =\ \lambda _{0} \lambda _{1} \lambda _{2}.  
\label{tau-real702} 
\end{eqnarray}
Further, by making use of (\ref{realC6}),   
we define the following functions on ${\mathcal U}_{H}$: 
\begin{eqnarray}
\kappa _{0} &=&  \mathrm{tr} (H)\, \lambda _{0} 
- \lambda _{0}^{2} + (n-1)\, c ,\nonumber\\ 
\kappa _{1} &=&  \mathrm{tr} (H)\, \lambda _{1} 
- \lambda _{1}^{2} + (n-1)\, c , \label{2021-10-18-a}\\
\kappa _{2} &=&  \mathrm{tr} (H)\, \lambda _{2} 
- \lambda _{2}^{2} + (n-1)\, c .\nonumber
\end{eqnarray}
We note that  
$\kappa _{0}$, $\kappa _{1}$  and $\kappa _{2}$ 
are eigenvalues of the Ricci operator $\mathcal{S}$ of $M$.
Now, using (\ref{real705}), (\ref{tau-real702})  
and (\ref{2021-10-18-a}) we can easy check that the following 
equations are satisfied: 
\begin{eqnarray}
\kappa _{0} - \kappa _{1}
&=&
( \lambda _{0} - \lambda _{1}) \, 
\left( \lambda _{2} - ( \phi -  \mathrm{tr} (H) ) \right) , 
\nonumber\\
\kappa _{0} - \kappa _{2}
&=&
( \lambda _{0} - \lambda _{2}) \, 
\left( \lambda _{1} - ( \phi -  \mathrm{tr} (H) ) \right) ,
\label{2022-01-17-a}\\
\kappa _{1} - \kappa _{2}
&=&
( \lambda _{1} - \lambda _{2}) \, 
\left( \lambda _{0} - ( \phi -  \mathrm{tr} (H) ) \right) , \nonumber
\end{eqnarray}
\begin{eqnarray}
& &
\left( \lambda _{0} - ( \phi -  \mathrm{tr} (H) ) \right)
\left( \lambda _{1} - ( \phi -  \mathrm{tr} (H) ) \right)
\left( \lambda _{2} - ( \phi -  \mathrm{tr} (H) ) \right)\nonumber\\
&=& 
\rho + \psi \, ( \phi - \mathrm{tr} (H)) 
+ \phi \, ( \phi - \mathrm{tr} (H))^{2} 
- ( \phi - \mathrm{tr} (H))^{3}\label{2022-01-17-b}\\  
&=& 
\rho + ( \phi - \mathrm{tr} (H)) \, 
( \psi +  \mathrm{tr} (H) ( \phi - \mathrm{tr} (H)) ) 
\ =\ \tau .\nonumber 
\end{eqnarray}
Now (\ref{2022-01-17-a}) and (\ref{2022-01-17-b}) lead to
\begin{eqnarray}
( \kappa _{0} - \kappa _{1})
  ( \kappa _{0} - \kappa _{2}) 
\, ( \kappa _{1} - \kappa _{2})  
&=&
(\lambda _{0} - \lambda _{1}) \, ( \lambda _{0} - \lambda _{2}) 
\, ( \lambda _{1} - \lambda _{2}) \, \tau .
\label{2022-01-17-c}  
\end{eqnarray}	
Thus we have														
\begin{theorem} 
Let $M$ be a hypersurface in a Riemannian space of constant curvature 
$N^{n+1}(c)$, $n \geq 3$,
having at every point $x \in {\mathcal U}_{H} \subset M$ 
three distinct principal curvatures  
$\lambda _{0}$, $\lambda _{1}$ and $\lambda _{2}$,
with multiplicities $n_{0}$, $n_{1}$ and $n_{2}$, respectively. 
Moreover, let the eigenvalues
$\kappa _{0}$, $\kappa _{1}$  and $\kappa _{2}$  
of the Ricci operator $\mathcal{S}$ of $M$
satisfy (\ref{2021-10-18-a}).
Then the function $\tau$, defined by (\ref{real705}), satisfies 
(\ref{2021-11-15-a})
on 
${\mathcal U}_{H}$. 
\end{theorem} 

\begin{remark}
Let $M$ be a hypersurface in a Riemannian space of constant curvature
$N^{n+1}(c)$, $n \geq 3$, satisfying assumptions of the last theorem.
\newline
(i) 
If we would have $\kappa _{0} = \kappa _{1} = \kappa _{2}$ at 
a point $x \in {\mathcal U}_{H}$ then $S = (\kappa / n)\, g$ at $x$,
a contradiction.
\newline
(ii) If $\kappa _{0} = \kappa _{1} \neq \kappa _{2}$,
or $\kappa _{0} = \kappa _{2} \neq \kappa _{1}$,
or $\kappa _{1} = \kappa _{2} \neq \kappa _{0}$, 
at a point $x \in {\mathcal U}_{H}$
then $\tau = 0$ and $M$ is partially Einstein at $x$.
\newline
(iii) If $\kappa _{0} \neq \kappa _{1}$ 
and $\kappa _{0} \neq \kappa _{2}$
and $\kappa _{1} \neq \kappa _{2}$ at a point $x \in {\mathcal U}_{H}$
then $\tau \neq 0$ at every point of a some neighbourhood 
$U \subset {\mathcal U}_{H}$ of $x$.
\end{remark}


\begin{example} (i)
Let $M$ be an austere hypersurface in a Riemannian space of constant curvature
$N^{n+1}(c)$, $n \geq 3$, with three distinct principal curvatures 
$\lambda _{0}$, $\lambda _{1}$ and $\lambda _{2}$.
It is known 
that at every point of such hypersurface there are exactly three
distinct principal curvatures
$\lambda _{0} = 0$, $\lambda _{1} = \lambda$ 
and $\lambda _{2} = - \lambda$,
with multiplicities $n - 2p$, $p$ and $p$, respectively,
where $1 \leq p < n$ and $\lambda$ is a function on $M$ 
(see, e.g., \cite{Bryant-1991}).
From  (\ref{2021-11-15-a})  and (\ref{2021-10-18-a})   
we get immedia\-tely $\kappa _{0} = (n-1)c$
and $\kappa _{1} = \kappa _{2} = (n-1)c - \lambda^{2}$, 
$\kappa _{0} - \kappa _{1} = \lambda^{2} >0$, 
and in a consequence $\tau = 0$
at every point of $M$. 
Now, by making use of (\ref{real702}), (\ref{realC6bb}) and  
(\ref{real703})--(\ref{2021.11.30.a5}), we obtain
\begin{eqnarray*}
\phi &=&  \rho = 0,\ \ \ \psi \ =\ \lambda ^{2},\ \ \
\mathrm{tr} (H)\ =\ 0,\ \ \ 
\mathrm{tr} (H^{2})\ =\ 2p \lambda ^{2},\ \ \
H^{3} \ =\ \lambda ^{2} H,\\  
\beta _{1} &=& \beta _{3}\ =\   \beta _{4}\ =\ 0,\ \ \
\beta_{2}\ =\ \frac{n - 2p - 1}{(n-2)(n-1)}\, \lambda ^{2},\\ 
\beta_{5} &=&  (n-2) \beta_{2} - (n-3)(n-1) c,\ \ \
\beta_{6}\ =\ \beta_{2} - (n-3) c , \\
\kappa &=& n (n-1) c - 2p \lambda ^{2},\ \ \
\mu \ =\ - \frac{2p}{(n-2)(n-1)}\, \lambda ^{2},\\
S^{2} &=& \left(2(n-1)c - \lambda^{2} \right) S 
+ (n-1)c \left( \lambda^{2} - (n-1)c  \right) g ,\\
R \cdot S &=&  c\, Q(g,S).
\end{eqnarray*} 
Thus $M$ is a Ricci-pseudosymmetric partially Einstein hypersurface.
Moreover,  the following conditions are satisfied on $M$, provided that $n \geq 4$,
\begin{eqnarray*}
C \cdot S &=&  \beta _{2}\, Q(g,S),\\
(n-2)\, R \cdot C  &=&   (n-2)\, Q(S,R) - (n-2)^{2}c\, Q(g,R)
- \frac{(n-3)c}{2} \, Q(S, g \wedge g) ,\\
(n-2)\, C \cdot R &=&
(n-3)\, Q(S,R) - \frac{(n - 3)c}{2} \, Q(S, g \wedge g)
\nonumber\\
& &
+  \left( \frac{\kappa }{n-1} +  \lambda ^{2}
-(n^{2} - 3n + 3) c \right) Q(g,R) ,\\
(n-2)\, C \cdot C &=&
(n-3)\, Q(S,R) 
+ \beta _{5}\, Q(g,R) + \frac{\beta _{6}}{2} Q(S, g \wedge g ) .
\end{eqnarray*}
The derivation-commutator $R \cdot C - C \cdot R$ is a linear combination
of the Tachibana tensors $Q(S,R)$ and $Q(g,R)$. Precisely, we have on $M$
\begin{eqnarray*}
(n-2)\, ( R \cdot C - C \cdot R) &=& Q(S,R)
+  \left( (n-1) c  -  \lambda ^{2} - \frac{\kappa }{n-1} \right) Q(g,R) \\
&=& Q\left(S + \left( (n-1) c  -  \lambda ^{2} - \frac{\kappa }{n-1} \right) g, 
R \right)  .
\end{eqnarray*}
We refer to \cite{DGHP-TZ 2} for a survey on manifolds (hypersurfaces)
with the derivation-commutator $R \cdot C - C \cdot R$ expressed as a linear combination
of Tachibana tensors.
\newline
(ii) (cf. {\cite[Section 10.2] {DHV2008}})
If $p = 1$ then, in view of results of \cite{1995_DDDVY},
(\ref{pseudo}) with $L_{R} = c$ 
holds on $M$.
Thus $M$ is a pseudosymmetric manifold. 
In particular, 
generalized Cartan hypersurfaces \cite{CHEN-1994} 
as well as 
the $3$-dimensional Cartan hypersurfaces are pseudosymmetric
(see also {\cite[Example 1] {DVY-1994}}).
If $p \geq 2$ then  $R \cdot S = c\, Q(g,S)$ holds on $M$,
i.e., $M$ is a non-pseudosymmetric Ricci-pseudosymmetric manifold.
In particular, Cartan hypersurfaces with 
$p = 2$ and $n = 6$, or
$p = 4$ and $n = 12$, or
$p = 8$ and $n = 24$,
are such manifolds. 
Pseudosymmetry type curvature properties of Cartan hypersurfaces
are presented in {\cite[Theorem 4.3] {DG90}}.
\newline
(iii) In {\cite[Example 3] {2005-LSchS}}
an example of a hypersurface $M$ in a $5$-dimensional unit 
sphere with four distinct principal curvatures is presented.
We can easy verify that $M$ is a partially Einstein hypersurface. 
\end{example}

Let (\ref{pseudo}) be satisfied on the set $\mathcal{U}_{R}$
of a hypersurface $M$ in a Riemannian space of constant curvature $N^{n+1}(c)$, 
$n \geq 4$. From {\cite[Theorem 5.1] {P39}} it follows that 
at every point $x \in \mathcal{U}_{R}$ we have: 
(i) 
$H^{2} = \alpha \, H + \beta \, g$, for some 
$\alpha , \beta \in \mathbb{R}$, 
i.e., 
$M$ has exactly two distinct principal curvatures 
$\lambda _{0}$ and $\lambda _{1}$,
or, 
(ii) 
$\mathrm{rank} (H) = 2$,
i.e., 
$M$ has principal curvatures 
$\lambda _{0} \neq 0$, $\lambda _{1} \neq 0$
and $\lambda _{2} = \lambda _{3}  =  \ldots = \lambda _{n-1} = 0$.

\begin{example} 
Let $M$ be a hypersurface in a Riemannian space of constant curvature 
$N^{n+1}(c)$, $n \geq 4$, having at every point principal curvatures
$\lambda _{0} \neq 0$, $\lambda _{1} \neq 0$
and
$\lambda _{2} = \lambda _{3} = \ldots = \lambda _{n-1} = 0$.
At every point of $M$ we have $\mathrm{rank} (H) = 2$. Thus 
$M$ is a type number two hypersurface 
(see, e.g., \cite{CHEN-YILDIRIM-2015}).
In view of {\cite[Theorem 4.2] {1995_DDDVY}},
we also have 
$R \cdot R = c\, Q(g,R)$
on $M$.
\newline
(i) Let $\lambda _{0} = \lambda _{1}$ at a point $x \in M$. Using now
suitable formulas of Section 3 we get immediately
$H^{2} = \lambda _{0} \, H$ 
($\alpha = \lambda _{0}$ and $\beta = 0$),
$\mathrm{tr} (H) - \alpha 
= 2 \lambda _{0} - \lambda _{0} = \lambda _{0} \neq 0$, 
\begin{eqnarray*}
\kappa _{0} &=&  \mathrm{tr} (H)\, \lambda _{0} 
- \lambda _{0}^{2} + (n-1)\, c 
\ = \ \lambda _{0}^{2} + (n-1)c , \\ 
\kappa _{1} &=&  \mathrm{tr} (H)\, \lambda _{1} 
- \lambda _{1}^{2} + (n-1)\, c 
\ = \ \lambda _{0}^{2} + (n-1)c ,\\
\kappa _{2} &=&  \mathrm{tr} (H)\, \lambda _{2} 
- \lambda _{2}^{2} + (n-1)\, c 
\ = \ (n-1)c ,\\
H 
&=& 
\lambda _{0}^{-1} ( S - (n-1) c \, g) ,\\ 
R &=& \frac{\lambda _{0}^{-2}}{2 }
( S - (n-1) c \, g) \wedge ( S - (n-1) c \, g) + \frac{c}{2} \, g \wedge g ,\\
S^{2} 
&=&
\lambda _{0} ( \lambda _{0}^{2} + 2 (n-1)c)\, H + (n-1)^{2}c^{2} \, g\\
&=&
( \lambda _{0}^{2} + 2 (n-1)c)\, S
- 
(n-1)c ( \lambda _{0}^{2} +  (n-1)c)\, g .
\end{eqnarray*}    
We assume now that $\lambda _{0} = \lambda _{1}$ at every point of $M$.
Thus $M$ is a $2$-quasi-Einstein and partially Einstein hypersurface.
Moreover, 
in view of Proposition 3.2, the hypersurface $M$ is a Roter space. 
Precisely, we have on $M$
\begin{eqnarray}
R  &=& \frac{\phi }{2} \, ( S - (n-1) c \, g) \wedge ( S - (n-1) c \, g) +
\frac{c}{2}\, g \wedge g,\ \ \  \phi \ =\ \lambda _{0}^{-2}.
\label{2022.09.21.aa}  
\end{eqnarray}
(ii)  Let $\lambda _{0} \neq \lambda _{1}$ at a point $x \in M$. 
We note that the case $\lambda _{1} = - \lambda _{0}$
was considered in Example 5.4 (i) (the case $p = 1$).
Using 
suitable formulas of Sections 4 and 5 we get immediately:
\begin{eqnarray*}
& &
\mathrm{tr} (H)  = \lambda _{0} + \lambda _{1} ,\ \ \
\phi = \lambda _{0} + \lambda _{1},\ \ \   \psi = - \lambda _{0} \lambda _{1},\ \ \
\rho = 0, \\
& &
H^{3} = \mathrm{tr} (H)\, H^{2} + \psi \, H,\ \ \  \tau = 0,\ \ \ A = \tau \, H = 0,
\end{eqnarray*}
\begin{eqnarray*}
S^{2} = \rho _{1} \, S + \rho _{0}\, g,\ \ \ 
\rho _{0} = - (n-1)c \, (  \lambda _{0} \lambda _{1} + (n-1)^{2}c^{2} ),\ \ \
\rho _{1} =  \lambda _{0} \lambda _{1} + 2 (n-1) c ,
\end{eqnarray*}
and
\begin{eqnarray*}
\kappa _{0} &=&  \mathrm{tr} (H)\, \lambda _{0} 
- \lambda _{0}^{2} + (n-1)\, c 
\ = \ \lambda _{0} \lambda _{1} + (n-1)c , \\ 
\kappa _{1} &=&  \mathrm{tr} (H)\, \lambda _{1} 
- \lambda _{1}^{2} + (n-1)\, c 
\ = \ \lambda _{0} \lambda _{1} + (n-1)c ,\\
\kappa _{2} &=&  \mathrm{tr} (H)\, \lambda _{2} 
- \lambda _{2}^{2} + (n-1)\, c 
\ = \ (n-1)c .
\end{eqnarray*}  
We assume that $\lambda _{0} \neq \lambda _{1}$ at every point of $M$.
We note that $\mathrm{rank} (S - (n-1)c\, g) = 2$ 
at every point of $M$.
Thus $M$ is a $2$-quasi-Einstein and partially Einstein hypersurface.
Furthermore, in view of {\cite[Theorem 3.2 (iii)] {2018_DGZ}},
(\ref{2022.09.21.aa}) holds on $M$,
for some function $\phi$. 
As it was stated in {\cite[Remark 5.2] {P39}}),
if 
$\lambda _{1} = - \lambda _{0} \neq 0$
and $\lambda _{2} = \lambda _{3} = \ldots = \lambda _{n-1} = 0$
at every point of a hypersurface $M$ in $N^{n+1}(c)$, $n \geq 4$,
then it is a minimal hypersurface realizing Chen's equality 
{\cite[Lemma 3.2] {CHEN-1993}}.     
\end{example}

\section{Hypersufaces satisfying ${\mathcal A}^{3} = \phi {\mathcal A}^{2} + \psi {\mathcal A} + \rho I\!d$. The non-partially Einstein case.}

Let $M$ be a hypersurface in $N_{s}^{n+1}(c)$, $n \geq 4$,
with the tensor $H$ satisfying (\ref{real702}) 
on ${\mathcal U}_{H} \subset M$.
From (\ref{real70456}) and (\ref{real707}) it follows that
if the function $\tau$,
defined by (\ref{real705}),
vanishes at a point $x \in {\mathcal U}_{H}$ 
then $M$ is partially Einstein at $x$.
We assume now that $\tau $ is non-zero at a point 
$x \in {\mathcal U}_{H}$. 
From
(\ref{real70456})--(\ref{2022.01.11.a1}), (\ref{2022.01.10.99})
and (\ref{2022.01.10.j}) 
it follows that 
\begin{eqnarray}
H &=& \tau ^{-1} \, A 
\ =\ \tau ^{-1} \, (S^{2} - \rho_{1} S - \rho _{0} g) ,
\label{new04-real7088} \\
R 
&=& \frac{\varepsilon}{2}\tau ^{-2} \, A \wedge A + \frac{c}{2}\, g \wedge g \nonumber\\
&=& \frac{\varepsilon}{2} \tau ^{-2} \,
(S^{2} - \rho_{1} S - \rho _{0} g) 
\wedge (S^{2} - \rho_{1} S - \rho _{0} g) + \frac{c}{2}\, g \wedge g ,
\label{real706.dd}\\
R \cdot S
&=&
- \tau ^{-1} \beta _{1} \, Q(S, S^{2}) 
+ \tau ^{-1} ( \rho + (n-1) c \beta_{1} )\, Q(g, S^{2})\nonumber\\
& &
+ ( c - \tau ^{-1} \rho  \rho _{1} 
- \tau ^{-1} \beta _{1} ( \rho _{0} + (n-1) c  \rho _{1}) )\, Q(g,S) ,
\label{2022.01.10.99new}\\
C \cdot S
&=&
- \tau ^{-1} \beta _{1}\, Q(S,S^{2}) 
+ \tau ^{-1} \beta _{4} \, Q(g,S^{2})\nonumber\\
& &
+
(\beta _{2} - \tau ^{-1} ( \beta _{1} \rho _{0} + \beta _{4} \rho _{1} ) )
\, Q(g,S) 
\label{2022.01.10.j.dd}
\end{eqnarray}
on some neighbourhood $U \subset {\mathcal U}_{H}$ of $x$.  
Furthermore, (\ref{realC6}) and (\ref{new04-real7088}) give
\begin{eqnarray}
H^{2} &=& \mathrm{tr}(H)\, H 
+ \varepsilon (n-1) c\, g - \varepsilon\, S \nonumber\\
&=& 
\tau^{-1} \mathrm{tr}(H)\, (S^{2} - \rho _{1}\, S - \rho _{0}\, g) 
+ \varepsilon (n-1) c\, g - \varepsilon \, S  \nonumber\\
&=& \tau^{-1} \mathrm{tr}(H) \, S^{2} 
- ( \tau^{-1} \mathrm{tr}(H) \rho _{1} + \varepsilon )\, S 
+ (\varepsilon (n-1) c - \tau^{-1} \rho _{0} \mathrm{tr}(H) )\, 
g \nonumber\\
&=& 
\tau^{-1} \left(
 \mathrm{tr}(H) S^{2} 
- ( \mathrm{tr}(H) \rho _{1} + \varepsilon \tau ) S 
+ (\varepsilon (n-1) c \tau  - \rho _{0} \mathrm{tr}(H) )\, g \right) .
\label{new04-real7099} 
\end{eqnarray}
Now
(\ref{new04-real7077}), 
by making use of (\ref{new04-real7088}) and (\ref{new04-real7099}),
turns into 
$S^{3} = \tau _{2}\, S^{2} + \tau _{1}\, S + \tau _{0}\, g$,
where $\tau _{0}, \tau _{1}, \tau _{2}$
are some functions on $U$. 
Thus we have
\begin{theorem}
Let $M$ be a hypersurface in a semi-Riemannian space of constant curvature
$N_{s}^{n+1}(c)$, $n \geq 4$,
with the tensor $H$ satisfying on ${\mathcal U}_{H}$
equation (\ref{real702})
and let  $\tau$ be the function defined on ${\mathcal U}_{H}$
by (\ref{real705}). We have:
\newline
(i) if $\tau$ vanishes at a point 
$x \in {\mathcal U}_{H}$ 
then $M$ is partially Einstein at $x$,
\newline
(ii) 
if $\tau$ is non-zero at a point 
$x \in {\mathcal U}_{H}$ then 
(\ref{new04-real7088})--(\ref{2022.01.10.j.dd})
hold on some
neighbourhood $U \subset {\mathcal U}_{H}$ of $x$
and at every point of $U$ the tensor $S^{3}$
is a linear combination of the tensors  
$g$, $S$ and $S^{2}$.
\end{theorem}

\begin{example}
(i) 
Let $M$ be a hypersurface in a Riemannian space of constant curvature
$N^{n+1}(c)$, $n \geq 3$,
having at every point $x \in {\mathcal U}_{H} \subset M$
three distinct principal curvatures 
$\lambda _{0}$, $\lambda _{1}$ 
and $\lambda _{2}$,
with multiplicities 
$1$, $1$ and $(n-2)$,
respectively, i.e.,
$M$ is $2$-quasi-umbilical at every point 
of ${\mathcal U}_{H}$ 
satisfying 
(\ref{real702})
(see, e.g., \cite{{CHEN-VERSTRAELEN-1976}, {DeVerYap}, {V2}}).
Further, from (\ref{2022-01-17-a}) and (\ref{2022-01-17-c}) 
we get immediately
\begin{eqnarray}
\tau &=&
(n-2) ( \lambda _{0} + (n-3) \lambda _{2})
 ( \lambda _{1} + (n-3) \lambda _{2})  \lambda _{2} .
\label{2022.02.04.a} 
\end{eqnarray}
(ii) Curvature properties of pseudosymmetry type of some particular
$2$-quasi-umbilical hypersurfaces $M$ in an Euclidean space
$\mathbb{E}^{n+1}$, $n \geq 4$, are given in
{\cite[Section 4] {DGP-TV02}}. 
At every point of the set ${\mathcal U}_{H}$
of those hypersurfaces
principal curvatures are the following
$\lambda _{0} = 0$, $\lambda _{1} = - (n-2) \lambda$ 
and $\lambda _{2} = \lambda$, where $\lambda \neq 0$. 
Therefore (\ref{2022.02.04.a}) turns into
$\tau = - (n-3) (n-2) \lambda ^{3}$.
Thus $M$ is a non-partially Einstein hypersurface. 
\newline 
(iii)
Let $M$ be a hypersurface in a Riemannian space of constant curvature
$N^{n+1}(c)$, $n \geq 4$,
having at every point $x \in {\mathcal U}_{H} \subset M$ 
three distinct principal curvatures  
$\lambda _{0} = a$, $\lambda _{1} = b$ and $\lambda _{2} = a+b$,
with multiplicities $1$, $1$ and $n - 2$, respectively.
Now the function $\tau$ takes the form   
\begin{eqnarray*}
\tau &=&
(n-2) ((n-2) a + (n-3) b)
 ( (n-3) a +  (n-2) b )  (a+b) .
\end{eqnarray*}
Evidently, if
$a \neq - b$ and 
$(n-2) a \neq - (n-3) b$
and
$(n-3) a \neq -  (n-2) b$
at a point $x \in {\mathcal U}_{H}$ 
then $\tau \neq 0$ at $x$
and 
$M$ is non-partially Einstein at $x$.
\end{example}

Let $M$ be a hypersurface in a Riemannian space of constant curvature
$N^{n+1}(c)$, $n \geq 4$,
having at every point $x \in {\mathcal U}_{H} \subset M$ 
three distinct principal curvatures  
$\lambda _{0}$, $\lambda _{1}$ and $\lambda _{2}$,
with multiplicities $1$, $n_{1}$ and $n_{2}$, respectively. 
We can check that 
(\ref{real705}) and (\ref{tau-real702}) yield 
(cf. {\cite[Section 7, eq. (7.22)] {2021-DGH}})
\begin{eqnarray}
\tau &=& 
\left( \lambda _{0}
 + (n_{1} - 1) \lambda _{1} + (n_{2} - 1) \lambda _{2} \right) 
\cdot \nonumber\\
& &\cdot \left( 
\lambda _{1} \lambda_{2} 
+ (n_{1} \lambda _{1} + n_{2} \lambda _{2})
((n_{1} - 1) \lambda _{1} + (n_{2} - 1) \lambda _{2} ) 
\right) .
\label{tau-real705} 
\end{eqnarray}

\begin{example} (i)
Hypersurfaces isometrically immersed 
in ${\mathbb{E}}^{n+1}$, $n \geq 5$, 
having at every point three distinct principal curvatures $\lambda _{0}$, 
$\lambda _{1}$ and $\lambda _{2}$, 
with multiplicities $1$, $n_{1}$ and $n_{2}$, respectively, 
were investigated for instance in 
\cite{{ABPRS}, {Sato-2002}, {Sato-Neto-2006}},
{see also 
\cite{FYZ} and 
\cite[Chapter 5.2] {OU-CHEN-2020}} and references therein).
\newline
(ii) Curvature properties of pseudosymmetry type 
of hypersurfaces in ${\mathbb{E}}^{n+1}$, $n \geq 5$,
having at every point three distinct principal curvatures $\lambda _{0}$, 
$\lambda _{1}$ and $\lambda _{2}$, 
with multiplicities $1$, $n_{1} = n_{2} \geq 2$,
are given
in {\cite[Section 4] {Saw-2015}}.
Moreover, in {\cite[Section 5, Example 5.1] {Saw-2015}},
by an application of results of {\cite[Section 2] {Sato-2002}}
(see also {\cite[Section 2] {ABPRS}} or 
{\cite[Section 2] {Sato-Neto-2006}})
a class of hypersurfaces 
$M$ in ${\mathbb{E}}^{n+1}$, $n \geq 5$, 
having at every point 
three distinct non-zero principal curvatures $\lambda _{0}$, 
$\lambda _{1}$ and $\lambda _{2}$, 
with multiplicities $1$, $n_{1} \geq 2$, $n_{2} \geq 2$, respectively,
was determined. 
In that class of hypersurfaces there are hypersurfaces
satisfying:
$n_{1} \neq  n_{2}$,
$( n_{1} - 1 ) \lambda_{1} + ( n_{2} - 1 ) \lambda_{2} = 0$
and  
$\phi = \lambda _{0} + \lambda _{1} + \lambda _{2} = \mathrm{tr} (H)$
(see {\cite[Section 5, Example 5.1 (ii)] {Saw-2015}}).
In addition,  
(\ref{tau-real702}) 
and 
(\ref{tau-real705})
yield
$\rho = \tau = \lambda_{0} \lambda_{1} \lambda_{2}$.  
On certain hypersurfaces of that class of hypersurfaces we also have:
$n_{1} = n_{2}$,
$\lambda_{1} = - \lambda_{2}$ and 
$\phi = \lambda _{0}  = \mathrm{tr} (H)$
(see {\cite[Section 5, Example 5.1 (iii)] {Saw-2015}}).
In addition, 
(\ref{tau-real702}) 
and 
(\ref{tau-real705})
yield
$\rho = \tau = - \lambda_{0} \lambda_{1}^{2} \neq 0$.  
\end{example}

\noindent
{\bf{Acknowledgments.}} 
The first two authors of this paper are supported 
by the Wroc\l aw University of Environmental and Life Sciences (Poland).

\end{document}